\documentclass[12pt]{article}
\usepackage{amsfonts}
\usepackage{amsmath,amssymb}
\usepackage{mathrsfs}
\usepackage{color}
\usepackage{amsxtra}
\usepackage{amstext}
\usepackage{amscd}
\usepackage{latexsym}
\usepackage{dsfont}
\usepackage{bbold}

\newtheorem{thm}{Theorem}[section]

\newtheorem{lem}[thm]{Lemma}
\newtheorem{cor}[thm]{Corollary}
\newtheorem{remark}[thm]{Remark}
\newtheorem{definition}[thm]{Definition}
\numberwithin{equation}{section}

\def\Tr{{\rm Tr}}
\def\tr{{\rm tr}}

\def\det{{\rm det}}
\def\Bbb{\mathbb}
\def\R{\Bbb R}
\def\C{\Bbb C}
\def\Z{\Bbb Z}

\newcommand{\thedate}
\begin{document}
%\hfill \textbf{\red{version JMAA}}

\[
\mbox {\bf\large A Dynamics Driven by Repeated Harmonic Perturbations }
\]

\vspace{1cm}

\begin{center}

\setcounter{footnote}{0}
\renewcommand{\thefootnote}{\arabic{footnote}}

\textbf{Hiroshi Tamura} \footnote{tamurah@staff.kanazawa-u.ac.jp}\\
 {Institute of Science and Engineering}      \\
              and \\
 Graduate School of the Natural Science and Technology\\
 Kanazawa University,\\
 Kanazawa 920-1192, Japan

\bigskip

\textbf{Valentin A.Zagrebnov }\footnote{Valentin.Zagrebnov@univ-amu.fr}\\
Institut de Math\'{e}matiques de Marseille - UMR 7373 \\
CMI-AMU, Technop\^{o}le Ch\^{a}teau-Gombert\\
39, rue F. Joliot Curie, 13453 Marseille Cedex 13, France\\
and \\
D\'{e}partement de Math\'{e}matiques \\
Universit\'{e} d'Aix-Marseille - Luminy, Case 901\\
163 av.de Luminy, 13288 Marseille Cedex 09, France

\vspace{1.5cm}

ABSTRACT

\end{center}

%%%%%%%%%%%%%%%%%%%%%%%%%%%%%%%%%%%%%%%%%%%%%%%%%%%%%%%%%%%%%%%%%%%%%%%%%%%%%%
We propose an exactly soluble $W^*$-dynamical system generated by repeated
harmonic perturbations of the one-mode quantum oscillator.
In the present paper we deal with the case of isolated system. Although dynamics is Hamiltonian and
quasi-free, it produces relaxation of initial state of the system to the steady state in the large-time limit.
The relaxation is accompanied by the entropy production and we found explicitly the rate for it.
Besides, we study evolution of subsystems to elucidate their eventual correlations and convergence
to equilibrium state. Finally we prove a universality of the dynamics driven by repeated harmonic perturbations
in a certain short-time interaction limit.
%\red{We will discuss the case of open system in another paper.}
%%%%%%%%%%%%%%%%%%%%%%%%%%%%%%%%%%%%%%%%%%%%%%%%%%%%%%%%%%%%%%%%%%%%%%%%%%%%%%
\newpage
%%%%%%%%%%%%%%%%%%%%%%%%%%%%%%%%%%%%%%%%%%%%%%%%%%%%%%%%%%%%%%%%%%%%%%%%%%%%%%
\tableofcontents
%%%%%%%%%%%%%%%%%%%%%%%%%%%%%%%%%%%%%%%%%%%%%%%%%%%%%%%%%%%%%%%%%%%%%%%%%%%%%%
%%%%%%%%%%%%%%%%%%%%%%%%%%%%%%%%%%%%%%%%%%%%%%%%%%%%%%%%%%%%%%%%%%%%%%%%%%%%%%
\section{Preliminaries and the Model}
%%%%%%%%%%%%%%%%%%%%%%%%%%%%%%%%%%%%%%%%%%%%%%%%%%%%%%%%%%%%%%%%%%%%%%%%%%%%%%
\subsection{Setup}\label{setup}
%%%%%%%%%%%%%%%%%%%%%%%%%%%%%%%%%%%%%%%%%%%%%%%%%%%%%%%%%%%%%%%%%%%%%%%%%%%%%%
We consider quantum system (one-mode quantum oscillator $\mathcal{S}$), which is successively perturbed
by time-dependent identical repeated interactions. This sequence of perturbation is switched on at the moment $t=0$
and it acts constantly on the interval $0 \leq t <\infty$. It is a common fashion to present this sequence as repeated
interactions of the system $\mathcal{S}$ with an \textit{infinite} time-equidistant \textit{chain}:
$\mathcal{C} = \mathcal{S}_1 + \mathcal{S}_2 + \ldots $, of subsystems $\{\mathcal{S}_{k}\}_{k\geq 1}$.
This visualisation is also motivated by certain physical models \cite{BJM}.

Below we suppose that the states of $\mathcal{S}$ and of every $\mathcal{S}_{k}$ are \textit{normal}, i.e. defined by
the \textit{density matrices} $\rho_0$ and $\{\rho_k\}_{k=1}^{\infty}$ on the Hilbert spaces
$\mathscr{H}_{\mathcal{S}}$ and $\{\mathscr{H}_{\mathcal{\mathcal{S}}_{k}}\}_{k=1}^{\infty}$, respectively.
The Hilbert space of the total system is then the tensor product $\mathscr{H}_{\mathcal{S}} \otimes \mathscr{H}_{\mathcal{C}}$.
Here the infinite product $\mathscr{H}_{\mathcal{C}}= \otimes_{k\geq 1} \mathscr{H}_{\mathcal{S}_{k}}$ stays for the Hilbert
space chain.

Since for any fixed moment $t\geq 0$, only a \textit{finite} number $N(t)$ of repeated interactions are involved into the
dynamics, {the subsystems $\{\mathcal{S}_{k}\}_{k> N(t)}$ are still \textit{independent} for different $k$, as well as
they are independent of components $\mathcal{S}$ and $\{\mathcal{S}_{k}\}_{k=1}^{N(t)}$.
On the other hand, the problem of correlations between components $\mathcal{S}$, $\mathcal{S}_{k}$
for $k \leq N(t)$  and between $\mathcal{S}_{k}$, $\mathcal{S}_{k^\prime}$ for $1 \leq k < k^\prime \leq N(t)$
is considered in Section \ref{subsystem}.
This peculiarity of repeated interactions allows to reduce the analysis of dynamics to the
finite tensor product: $\mathscr{H}_{\mathcal{C}_N }= \otimes_{k=1}^{N}\mathscr{H}_{\mathcal{S}_{k}}$.
Then one recovers the above infinite chain $\mathcal{S} + \mathcal{C}$ as the limit $N \rightarrow \infty$,
\textit{a posteriori}}.

%%%%%%%%%%%%%%%%%%%%%%%%%%%%%%%%%%%%%%%%%%%%%%%%%%%%%%%%%%%%%%%%%%%%%%%%%%%%%%%%%%%%%%%%%%%%%%%%%%%%%%%%%%%%%%%%
Details of dynamics are presented in the next Section \ref{HDIS}. In this section we mention our
{\textit{guiding}} hypothesises.

\smallskip

\noindent
\textbf{Hypothesis 1:} For $t\leq 0$, \textit{all} components of $\mathcal{S}$ and
$\{\mathcal{S}_{k}\}_{k=1}^{N}$ are \textit{independent}, i.e.
the state of $\mathcal{S} + \mathcal{C}_N$ is described as a finite tensor product:
$\omega_{\mathcal{S} + \mathcal{C}_N}: = \omega_{\mathcal{S}} \otimes \bigotimes_{k=1}^N \omega_{\mathcal{S}_{k}}$.
We suppose that each of the state in the product is \textit{normal}.
%%%%%%%%%%%%%%%%%%%%%%%%%%%%%%%%%%%%%%%%%%%%%%%% Remark %%%%%%%%%%%%%%%%%%%%%%%%%%%%%%%%%%%%%%%%%%%%%%%%%%%%%%%%
\begin{remark}\label{Rem1}
Although it is {not} decisive for our arguments, we recall that the product $\mathscr{H}_{\mathcal{C}}$
as well as the von Neumann algebra of observables of the infinite total system
$\mathfrak{M}= \mathfrak{M}_{\mathcal{S}} \otimes \mathfrak{M}_{\mathcal{C}}$ can be correctly defined, see e.g. \cite{BR1}
(Sections 2.7.2 and 2.7.3). Here $\mathfrak{M}_{\mathcal{S}} \subseteq \mathcal{L}(\mathscr{H}_{\mathcal{S}})$ and
$\mathfrak{M}_{\mathcal{C}} = \otimes_{k\geq 1} \mathfrak{M}_{\mathcal{S}_{k}}\subseteq \mathcal{L}(\mathscr{H}_{\mathcal{C}})$
are von Neumann algebras of bounded operators $\mathcal{L}$ on  $\mathscr{H}_{\mathcal{S}}$ and $\mathscr{H}_{\mathcal{C}}$,
respectively.
\end{remark}
%%%%%%%%%%%%%%%%%%%%%%%%%%%%%%%%%%%%%%%%%%%%%%%%%%%%%%%%%%%%%%%%%%%%%%%%%%%%%%%%%%%%%%%%%%%%%%%%%%%%%%%%%%%%%%%%%

A basic ingredient in construction of dynamical system $\mathcal{S} + \mathcal{C}_N$ is the \textit{one-mode} quantum
harmonic oscillator. Recall that it can be described by (\textit{unbounded}) boson annihilation and creation operators $a, a^*$
defined in the Fock space $\mathscr{F}$. They realise a representation of the Canonical Commutation Relations
(CCR) in $\mathscr{F}$, i.e. formally satisfy the operator relations:
\[
      [a, a^*] = \mathbb{1}, \quad  [a, a] = 0,
        \quad [a^*, a^*] = 0 \ .
\]
Here $\mathbb{1}$ denotes the unit operator on $\mathscr{F}$.
Let $\Omega\in \mathscr{F}$ be the vacuum vector: $a \, \Omega = 0$. Then the Hilbert space space  $\mathscr{F}$ is
vector-norm {completion} of the algebraic span of vectors $\{(a^{*})^m \Omega\}_{m\geq 0}$.

Denote by $\{\mathscr{H}_k\}_{k=0}^{N}$, the copies of the Fock space $\mathscr{F}$ for \textit{arbitrary} but
finite $N \in \mathbb{N}$ and by $\mathscr{H}^{(N)}$, the Hilbert space
$ \mathscr{H}_{\mathcal{S}} \bigotimes \mathscr{H}_{\mathcal{C}_N} $ for $\mathcal{S} + \mathcal{C}_N $ with
$\mathscr{H}_{\mathcal{S}} = \mathscr{H}_{0}$ and $ \mathscr{H}_{\mathcal{S}_k} =  \mathscr{H}_{k} \; ( k= 1, \cdots, N)$, i.e.,
\begin{equation}\label{H-space}
\mathscr{H}^{(N)}:=
 \mathscr{H}_{0} \otimes \bigotimes_{k=1}^{N} \mathscr{H}_k = \mathscr{F}^{\otimes(N+1)} \ .
\end{equation}

Now we define in (\ref{H-space}) two operators
\begin{equation}\label{S-bosons}
b_0 := a\otimes \mathbb{1} \otimes  \ldots \otimes \mathbb{1}, \qquad
b_0^* := a^*\otimes \mathbb{1} \otimes  \ldots \otimes \mathbb{1}  \ .
\end{equation}
These CCR boson operators serve for description of the system $\mathcal{S}$ with Hamiltonian
\begin{equation}\label{H-Cav-S}
 H_{\mathcal{S}}:= E \, b_0^*b_0 \ \ , \ \ {\rm{dom}}(H_{\mathcal{S}}) \subset \mathscr{H}^{(N)} \ , \ E > 0 \ .
\end{equation}
It is a one-mode harmonic oscillator with discrete spectral spacing $E$.
{We consider it as an isolated (\textit{ideal}) one-mode quantum harmonic subsystem.}

The subsystems $\{\mathcal{S}_{k}\}_{k\geq 1}$ that we consider in the present paper are in turn identical one-mode harmonic
quantum oscillators with discrete spectral spacing $\epsilon$. Then to describe evolution of the finite system
$\mathcal{S} + \mathcal{C}_{N}$ due to $N = N(t)$ consecutive interactions, we define in space (\ref{H-space}) the
sequence of boson operators  $\{b_k , b_{k}^* \}_{k=1}^{N}$ in the space (\ref{H-space}):
\[
     b_1:= \mathbb{1} \otimes  a\otimes \mathbb{1} \otimes \mathbb{1} \otimes  \ldots \otimes \mathbb{1} ,
   \qquad
      b_1^* := \mathbb{1}\otimes a^*\otimes \mathbb{1} \otimes \mathbb{1} \otimes  \ldots \otimes \mathbb{1} ,
\]
\[
     b_2:= \mathbb{1}\otimes \mathbb{1} \otimes  a\otimes \mathbb{1} \otimes  \ldots \otimes \mathbb{1} ,
\qquad
      b_2^* := \mathbb{1}\otimes \mathbb{1}\otimes a^*\otimes \mathbb{1} \otimes  \ldots \otimes \mathbb{1} ,
\]
\begin{equation}\label{S-bosons1}
\ldots \ldots \ldots \ldots \ldots \ldots \ldots \ldots \ldots \ldots
\ldots \ldots \ldots \ldots \ldots \ldots \ldots \ldots \ldots
\end{equation}
\[
     b_N:= \mathbb{1}\otimes \mathbb{1} \otimes  \mathbb{1} \otimes \mathbb{1} \otimes  \ldots \otimes a  ,
\qquad
     b_N^* := \mathbb{1}\otimes \mathbb{1}\otimes \mathbb{1}\otimes \mathbb{1} \otimes  \ldots \otimes a^* .
\]
The Hamiltonian of each of subsystem $\mathcal{S}_k$ has the form
\begin{equation}\label{H-Sk}
H_{\mathcal{S}_k} : =  \epsilon \, b_k^*b_k  \ , \ \
{\rm{dom}}(H_{\mathcal{S}_k}) \subset \mathscr{H}^{(N)} \ , \ \epsilon > 0
\ , \
k=1,2,\ldots,N \ .
\end{equation}
Altogether the boson operators (\ref{S-bosons}) and (\ref{S-bosons1}) formally satisfy the CCR in the space (\ref{H-space}):
\begin{equation}\label{CCR}
[b_k, b^*_{k^\prime}] = \delta_{k,k^\prime} \mathbb{1} , \quad
[b_k, b_{k^\prime}] =  [b^*_k, b^*_{k^\prime}] = 0 \ , \ k,k^\prime =0, 1, 2, \ldots, N .
\end{equation}
%%%%%%%%%%%%%%%%%%%%%%%%%%%%%%%%%%%%%%%%%%%%%% Remark %%%%%%%%%%%%%%%%%%%%%%%%%%%%%%%%%%%%%%%%%%%%%%%%%%%%%%%%%
\begin{remark}\label{Rem-Phys}
Note that there is some physical interpretation \cite{NVZ} behind of this mathematical modelling.
For example, the system $\mathcal{C}_N $ can be identified with a chain of $N$ quantum particles (atoms or molecules) (\ref{H-Sk})
with \textit{harmonic} internal degrees of freedom interacting \textit{one-by-one} with $E$-one-mode quantum cavity (\ref{H-Cav-S}).
\end{remark}
%%%%%%%%%%%%%%%%%%%%%%%%%%%%%%%%%%%%%%%%%%%%%%%%%%%%%%%%%%%%%%%%%%%%%%%%%%%%%%%%%%%%%%%%%%%%%%%%%%%%%%%%%%%%%%%

\smallskip
%%%%%%%%%%%%%%%%%%%%%%%%%%%%%%%%%%%%%%%%%%%%%%%%%%%%% H2 %%%%%%%%%%%%%%%%%%%%%%%%%%%%%%%%%%%%%%%%%%%%%%%%%%%%%%%%%%%%%%
\noindent
\textbf{Hypothesis 2:} (\textit{Tuned interaction}) We
consider repeated perturbations in the tuned regime: for any
moment $t\geq 0$ exactly {\textit{one}} subsystem $\mathcal{S}_n$ is interacting
with the system  $\mathcal{S}$ during a fixed time $\tau >0$. Here $n = [t/\tau] + 1$, where $[x]$ denotes the
\textit{integer} part of $x\geq 0$.
%%%%%%%%%%%%%%%%%%%%%%%%%%%%%%%%%%%%%%%%%%%%%%%%%%%%%%%%%%%%%%%%%%%%%%%%%%%%%%%%%%%%%%%%%%%%%%%%%%%%%%%%%%%%%%%%%%%%%%%

\smallskip
%%%%%%%%%%%%%%%%%%%%%%%%%%%%%%%%%%%%%%%%%%%%%%%%%%%%%% H3 %%%%%%%%%%%%%%%%%%%%%%%%%%%%%%%%%%%%%%%%%%%%%%%%%%%%%%%%%%%%%%
\noindent
\textbf{Hypothesis 3:} The time-dependent repeated interaction is a piecewise \textit{constant}
operator, which is taken as the sum over $\ n \geq 1$ of the following  \textit{bilinear} forms in operators
(\ref{S-bosons}), (\ref{S-bosons1}):
\begin{equation}\label{Int}
K_n(t):=\chi_{[(n-1)\tau, n\tau)}(t) \,  \eta \  (b_0^*b_n + \ b_n^*b_0) \ , \quad  \ {\eta >0 }.
\end{equation}
Here $\chi_{\mathcal{I}}(x)$ is the characteristic function of the set $\mathcal{I}$.

%%%%%%%%%%%%%%%%%%%%%%%%%%%%%%%%%%%%%%%%%%%%%%%%%%%%%% 1.2 %%%%%%%%%%%%%%%%%%%%%%%%%%%%%%%%%%%%%%%%%%%%%%%%%%%%%%%%%%%%%
\subsection{The Model}\label{model}
%%%%%%%%%%%%%%%%%%%%%%%%%%%%%%%%%%%%%%%%%%%%%%%%%%%%%%%%%%%%%%%%%%%%%%%%%%%%%%%%%%%%%%%%%%%%%%%%%%%%%%%%%%%%%%%%%%%%%%%%
For any $N\geq 1$, the Hamiltonian $H_{N}(t)$ of  \textit{non-autonomous} system $\mathcal{S}+\mathcal{C}_{N}$
is defined in the space (\ref{H-space}) as the sum of all ingredients (\ref{H-Cav-S}), (\ref{H-Sk})
and (\ref{Int}). This sum is essentially self-adjoint operator in $\mathscr{H}^{(N)}$. Since it does not produce any confusion,
we denote its \textit{closure} by the same symbol:
\begin{align}\label{Ham-Model}
H_{N}(t):=& \ H_{\mathcal{S}}+\sum_{k=1}^{N}(H_{\mathcal{S}_k}+K_k(t))\\
=& \ E b_0^*b_0 + \epsilon \sum_{k=1}^{N} b_k^*b_k + \sum_{k=1}^{N}\chi_{[(k-1)\tau, k\tau)}(t)\,
\eta \, (b_0^*b_k +  \ b_k^*b_0) \, . \nonumber
\end{align}

By virtue of (\ref{Int}), (\ref{Ham-Model}) only $\mathcal{S}_n$ interacts with
$\mathcal{S}$ {for $t\in[(n-1)\tau, n\tau)$, $n\geq 1$}, i.e. the system $\mathcal{S}+\mathcal{C}_{N}$ is
\textit{autonomous} on this time-interval with self-adjoint Hamiltonian
\begin{equation}\label{Ham-n}
H_n(b) := E \, b_0^*b_0 + \epsilon\sum_{k=1}^{N}b_k^*b_k + \eta \, (b_0^*b_n +
 \ b_n^*b_0) \ , \ n \leq N \ .
\end{equation}
Here again operator $H_n$ denotes the closure of the algebraic sum of operators in the right-hand side of (\ref{Ham-n}).
{We shall consider only the case $ t < N\tau $.}

Note that CCR (\ref{CCR}) and definition of Hamiltonian (\ref{Ham-n}) yield
\begin{eqnarray} \label{H-n-commut}
&&[H_n, b_0] = -Eb_0 -\eta b_n \  ,  \ [H_n, b_j] = -\epsilon b_j-\delta_{jn} {\eta}b_0 \ , \\
&& \ \ [H_n, b_0^*] = Eb_0^* + {\eta} b_n^* \  ,  \ [H_n, b_j^*] = \epsilon b_j^*+\delta_{jn} {\eta}b_0^* \ , \nonumber
\end{eqnarray}
for any $1 \leq j \leq N$.

Moreover, since Hamiltonian (\ref{Ham-n}) is bilinear, there exists a \textit{canonical} (i.e. {CCR-}preserving )
linear transformation \cite{Ar0}
\begin{equation}\label{canon-transform1}
P_n : \{b_k\}_{k=0}^{N} \rightarrow \{c_k\}_{k=0}^{N} \ ,
\end{equation}
which diagonalises (\ref{Ham-n}):
\begin{equation}\label{canon-transform2}
\widetilde{H}_n(c) := \varepsilon_{0} \, c_0^* c_0 +  \varepsilon_{1} \, c_1^* c_1
+ \sum_{k=2}^{N} \varepsilon_k  c_k^*c_k \ ,
\end{equation}
{where}
\begin{eqnarray}
&&\varepsilon_{0} := \frac{1}{2}[(E+\epsilon) + \sqrt{(E-\epsilon)^2 + 4 \eta^2}] \ ,  \label{Sp-canon-transform0}\\
&&\varepsilon_{1} := \frac{1}{2}[(E+\epsilon) - \sqrt{(E-\epsilon)^2 + 4 \eta^2}] \ ,  \label{Sp-canon-transform1}
\end{eqnarray}
and $\varepsilon_2 = \ldots = \varepsilon_N = \epsilon $.
%%%%%%%%%%%%%%%%%%%%%%%%%%%%%%%%%%%%%%%%%%%%%%%%%%%%%%%%%%%%%%%%%%%%%%%%%%%%%%%%%%%%%%%%%%%%%%%%%%%%%%%%%%%%%%%%%%%%%%%%%%

\smallskip
%%%%%%%%%%%%%%%%%%%%%%%%%%%%%%%%%%%%%%%%%%%%%%%%%%%% H4 %%%%%%%%%%%%%%%%%%%%%%%%%%%%%%%%%%%%%%%%%%%%%%%%%%%%%%%%%%%%%%%%%%
\noindent
\textbf{Hypothesis 4:} By virtue of (\ref{Sp-canon-transform0}), (\ref{Sp-canon-transform1}) to keep Hamiltonian
(\ref{Ham-Model}) (or (\ref{Ham-n})) semi-bounded from \textit{below}, we must {impose the condition}
\begin{equation}\label{H4}
\eta^2 \leq E \, \epsilon \ .
\end{equation}

%At this stage few remarks concerning our model (\ref{Ham-Model}) are in order.
%%%%%%%%%%%%%%%%%%%%%%%%%%%%%%%%%%%%%%%%%%%%%%% Remark %%%%%%%%%%%%%%%%%%%%%%%%%%%%%%%%%%%%%%%%%%%%%%%%%%%%%%%%
\begin{remark}\label{Rem2}
{\rm{(i)}} The  non-autonomous system $\mathcal{S}+\mathcal{C}$ formally corresponds to Hamiltonian $H_{\infty}(t)$,
$t\in [0, \infty)$, i.e., {the case $N=\infty$ of} (\ref{Ham-Model}). Then for $t\in[(n-1)\tau, n\tau)$, $n\geq1$, and for
any $j\in \mathbb{N}$ one obtains:
\begin{eqnarray} \label{H-inf-commut}
&&[H_{\infty}(t), b_0] = -E b_0 -\eta b_n \  ,  \ [H_{\infty}(t), b_j] = -\epsilon b_j-\delta_{jn} {\eta}b_0 \ , \\
&& \ \ [H_{\infty}(t), b_0^*] = Eb_0^* + {\eta} b_n^* \  ,
\ [H_{\infty}(t), b_j^*] = \epsilon b_j^*+\delta_{jn} {\eta}b_0^* \ ,
\nonumber
\end{eqnarray}
{corresponding} to (\ref{H-n-commut}). These formal calculations can be justified in the framework
of the infinite tensor products and the von Neumann algebras mentioned above \cite{BR1}. Instead of that, in the present paper
we look on $N = \infty$ as \textit{a posteriori} inductive limit.

{\rm{(ii)}} Below. we study the non-autonomous system $\mathcal{S}+\mathcal{C}_{N}$ for arbitrary but
fixed $N\geq 1$ conditioned by $t\in [0,\tau N)$. Then (\ref{H-n-commut}) and (\ref{H-inf-commut}) coincides. Hence,
we keep {using} the notation  $\mathcal{S}+\mathcal{C}$ {in this case}. We restore the skipped subindex $N$ only if the
finiteness of the chain $\mathcal{C}$ is indispensable to stress.

{\rm{(iii)}} Note that in contrast to the case studied in {\rm{\cite{NVZ}}}, the subsystems {$\mathcal{S}_k$} are
not rigid and the interaction (\ref{Int}) is inelastic. Then, besides the energy/entropy exchange between
$\mathcal{S}$ and $\mathcal{C}$, the repeated perturbations may produce {{\rm{entanglement, or correlations,}}} of states
in the subsystems $\{\mathcal{S}_k\}_{k=1}^{n(t)} \ $, where $n(t) = [t/\tau] +1 $.
\end{remark}
%%%%%%%%%%%%%%%%%%%%%%%%%%%%%%%%%%%%%%%%%%%%%%%%%%%%%%%%%%%%%%%%%%%%%%%%%%%%%%%%%%%%%%%%%%%%%%%%%%%%%%%%%%%%%%%%

We conclude this section by a general lemma based on harmonic structure of the model {and by the final Hypothesis 5}.
{The first} yields an explicit form for commutators (\ref{H-n-commut}) and canonical transformation
(\ref{canon-transform1}), (\ref{canon-transform2}).
%%%%%%%%%%%%%%%%%%%%%%%%%%%%%%%%%%%%%%%%%%%%%%%%% Lemma %%%%%%%%%%%%%%%%%%%%%%%%%%%%%%%%%%%%%%%%%%%%%%%%%%%%%%%%
\begin{lem} \label{lem11}
For $j = 0, 1, 2, \ldots, N $ and $n = 1, 2, \ldots, N $, one gets
\begin{equation}\label{U-t-*}
     e^{it H_n} b_j e^{-it H_n} = \sum_{k=0}^N(U_{n}^*(t))_{jk}b_k,
     \quad e^{it H_n} b_j^* e^{-it H_n}
     = \sum_{k=0}^N\overline{(U_{n}^*(t))_{jk}}b_k^* ,
\end{equation}
\begin{equation}\label{U-t}
     e^{-it H_n} b_j e^{it H_n} = \sum_{k=0}^N(U_n(t))_{jk}b_k,
     \quad e^{-it H_n} b_j^* e^{it H_n}
     = \sum_{k=0}^N\overline{(U_n(t))_{jk}}b_k^* \ ,
\end{equation}
for $t \geq 0$. Here $U_n(t)$ and $V_n(t)$ are $(N+1)\times(N+1)$ matrices related by
$U_n(t) = e^{it\epsilon}V_n(t)$, where
\begin{equation}\label{V-n}
       (V_n(t))_{jk} := {\begin{cases}
        g(t)z(t)\,\delta_{k0} + g(t)w(t)\,\delta_{kn} & \quad (j = 0) \\
    g(t){w(t)}\,\delta_{k0} + g(t){z(-t)}\,\delta_{kn} &
        \quad (j = n) \\
                   \qquad \delta_{jk}    & \quad (\mbox{otherwise})
                    \end{cases}}
\end{equation}
with
\begin{equation}\label{g-w}
      g(t) := e^{it(E-\epsilon)/2} \ , \
           w(t) :=\frac{2i\eta}{\sqrt{(E-\epsilon)^2+4\eta^2}}
          \sin t\sqrt{\frac{(E-\epsilon)^2}{4}+\eta^2} \ ,
\end{equation}
\begin{equation}\label{z}
          z(t) := \cos t\sqrt{\frac{(E-\epsilon)^2}{4}+\eta^2}
          +\frac{i(E-\epsilon)}{\sqrt{(E-\epsilon)^2+4\eta^2}}
          \sin t\sqrt{\frac{(E-\epsilon)^2}{4}+\eta^2} \ .
         \end{equation}
\end{lem}
%%%%%%%%%%%%%%%%%%%%%%%%%%%%%%%%%%%%%%%%%%%%%%%%%%%%%%%%%%%%%%%%%%%%%%%%%%%%%%%%%%%%%%%%%%%%%%%%%%%%%%%%%%%%%%
%%%%%%%%%%%%%%%%%%%%%%%%%%%%%%%%%%%%%%%%%%%%%%%%%% Remark %%%%%%%%%%%%%%%%%%%%%%%%%%%%%%%%%%%%%%%%%%%%%%%%%%%%
\begin{remark}\label{Rem3}
Note that by definitions (\ref{g-w}) and (\ref{z}), we get $|z(t)|^2  + |w(t)|^2 =1 $, $z(-t) = \overline{z(t)}$ and
$w(t) = - \overline{w(t)}$. Therefore, the matrix
\[
    M(t) :=  \begin{pmatrix} z(t) & w(t) \\
                              &      \\
                 {w(t)} & {z(-t)} \end{pmatrix}
\]
is unitary. For $N=1$, one gets $M(t)= \overline{g(t)} V_1(t)$, see (\ref{V-n}). Moreover, (\ref{U-t-*}) and (\ref{U-t})
imply that $\{ V_n(t) \}_{t\in\R}$ and $\{U_n(t)\}_{t\in\R}$ are in fact one-parameter groups of $ (N+1) \times (N+1) $
unitary matrices.
\end{remark}
%%%%%%%%%%%%%%%%%%%%%%%%%%%%%%%%%%%%%%%%%%%%%%%%%%%%%%%%%%%%%%%%%%%%%%%%%%%%%%%%%%%%%%%%%%%%%%%%%%%%%%%%%%%%%%
{\sl Proof }(of Lemma \ref{lem11} ): Let $\{J_n\}_{n=1}^N$ and $\{X_n\}_{n=1}^N$ be $(N+1)\times (N+1)$ Hermitian matrices
given by
\begin{equation}\label{J-n}
      (J_n)_{jk} := \begin{cases}
                      1 & \quad ( j=k=0 \; \mbox{ or } \; j=k=n) \\
                      0 & \quad \mbox{otherwise}
                   \end{cases},
\end{equation}
\begin{equation}\label{X-n}
     (X_n)_{jk} := \begin{cases}
                      (E-\epsilon)/2 & \quad (j,k)=(0,0)  \\
                      -(E-\epsilon)/2 & \quad (j,k)=(n,n)  \\
                      \eta & \quad (j,k)=(0,n) \\
                      \eta & \quad (j,k)=(n,0) \\
                      0 & \quad \mbox{otherwise}
                   \end{cases} \, .
\end{equation}
We define the matrices
\begin{equation}\label{Y-n}
    Y_n: =\epsilon I + \frac{E-\epsilon}{2}J_n + X_n \  \ (n= 1, \ldots ,N) \ ,
\end{equation}
where $I$ is the $(N+1)\times (N+1)$ identity matrix. Then Hamiltonian (\ref{Ham-n}) takes the form
\begin{equation}\label{HY}
H_n={\sum_{j,k=0}^{N} (Y_n)_{jk}b_j^*b_k} \ .
\end{equation}
Since $Y_n$ is Hermitian, there exists a \textit{diagonal} matrix $\Lambda$ and  unitary mapping
$P_n: \mathbb{R}^{N+1}\rightarrow \mathbb{R}^{N+1}$,  implemented by transformation (\ref{canon-transform1}),
such that $ Y_{n} = P_{n}^{*} \Lambda P_n$ holds. Recall that after canonical transformation the matrix
{$\Lambda := \{\Lambda_{ij} \}_{i,j =0}^{N}= \{\delta_{ij} \, \varepsilon_{j}\}_{i,j =0}^{N}$} is universal and independent
of $n$, see (\ref{canon-transform2}). Then new operators (\ref{canon-transform1})
\begin{equation}\label{c-c}
c_j=\sum_{k=0}^N (P_n)_{jk} \ b_k, \quad
c_j^*=\sum_{k=0}^N \overline{(P_n)_{jk}} \ b_k^*
\quad \; (j=0,1, \ldots, N)
\end{equation}
satisfy CCR in the space $\mathscr{H}^{(N)}$  (\ref{H-space}) and  diagonalise (\ref{HY}):
\begin{equation}\label{H-c-c}
\widetilde{H}_n = {\sum_{j=0}^{N}\Lambda_{jj}c_j^*c_j} \ ,
\end{equation}
where $ \Lambda_{jj} = \varepsilon_{j} $, see (\ref{canon-transform2}). Therefore, the set of all eigenvectors of
$\widetilde{H}_n$ is:
\begin{equation}\label{eigenvectorsH-n}
         \Big\{ \, \prod_{j=0}^N\frac{(c_j^{*})^{n_j}}{\sqrt{n_j!}}
         \, \Omega\otimes \ldots\otimes\Omega \,
          \Big| \, n_j\in \Z_+ \quad (j=0,1,\ldots, N) \, \Big\} \ .
\end{equation}
Note that it forms a complete orthonormal basis in $\mathscr{H}^{(N)}$.
The linear envelope $\mathscr{H}^{(N)}_0$ of the set (\ref{eigenvectorsH-n}) is  invariant subspace
for transformations $e^{it\widetilde{H}_n}$ and its norm-closure coincides with $\mathscr{H}^{(N)}$. Then by
(\ref{H-c-c}) one gets on (\ref{eigenvectorsH-n})
\[
         e^{it\widetilde{H}_n}c_je^{-it\widetilde{H}_n} = e^{-it\Lambda_{jj}}c_j, \quad
         e^{it\widetilde{H}_n}c_j^*e^{-it\widetilde{H}_n} = e^{it\Lambda_{jj}}c_j^* \ .
\]
Now taking into account canonical transformation (\ref{c-c}), we obtain
\[
    e^{itH_n} b_j e^{-itH_n} = \sum_{k=0}^N  (P_{n}^{*})_{jk} \, e^{it\widetilde{H}_n} c_k e^{-it\widetilde{H}_n}
\]
\begin{equation}\label{Y-t}
       = \sum_{k,l=0}^N {(P_n^*)}_{jk}e^{-it\Lambda_{kk}}(P_n)_{kl}b_l
       = \sum_{l=0}^N\big( e^{-itP_n^*\Lambda P_n} \big)_{jl}b_l
       = \sum_{l=0}^N\big( e^{-itY_n} \big)_{jl}b_l  \ .
\end{equation}
Similarly we get on $\mathscr{H}^{(N)}_0$:
\[
      e^{itH_n}b_j^*e^{-itH_n} =
    \sum_{l=0}^N\overline{\big( e^{-itY_n }\big)_{jl}}b_l^*  \ .
\]

Note that by virtue of (\ref{J-n}), (\ref{X-n}), one has identities
\[
        X_n^2 = \Big(\frac{(E-\epsilon)^2}{4} + \eta^2\Big) \, J_n
       \quad \mbox{and}\quad
         J_nX_n =X_n \ .
\]
Together with definition (\ref{Y-n}), they yield
\begin{equation}\label{Y-U}
          e^{itY_n} = e^{it\epsilon}
          \bigg( I - J_n + e^{it(E-\epsilon)/2}\Big\{
          J_n \cos \, t \, \sqrt{\frac{(E-\epsilon)^2}{4} + \eta^2}
\end{equation}
\begin{equation*}
         + iX_n \left[\frac{(E-\epsilon)^2}{4} + \eta^2 \right]^{-1/2}
          \sin \, t \, \sqrt{\frac{(E-\epsilon)^2}{4} + \eta^2} \;
           \Big\}\bigg)
          = e^{it\epsilon}V_n(t) = U_n(t) .
\end{equation*}
Inserting now (\ref{Y-U}) into (\ref{Y-t}), we prove (\ref{U-t-*}).

Since $ U_n(t)^* = U_n(-t)$, one can similarly establish (\ref{U-t}). \hfill$\square$

\smallskip

%%%%%%%%%%%%%%%%%%%%%%%%%%%%%%%%%%%%%%%%% H5%%%%%%%%%%%% %%%%%%%%%%%%%%%%%%%%%%%%%%%%
\noindent
{\textbf{Hypothesis 5:} The parameters of the model (\ref{Ham-Model}) and the interaction time $\tau$ verify in addition to
(\ref{H4}) the conditions $|w(\tau)|<1$ and  $|z(\tau)|<1$, which are satisfied for, e.g., small $\tau$:
\begin{equation}\label{H5}
\tau \sqrt{{(E-\epsilon)^2}/{4}+\eta^2} < \pi/2 ,
\end{equation}
see (\ref{g-w}), (\ref{z}).}
%%%%%%%%%%%%%%%%%%%%%%%%%%%%%%%%%%%%%%%%%%%%%%%%%% Remark %%%%%%%%%%%%%%%%%%%%%%%%%%%%%%%%%%%%%%%%%%%%%%%%%%%%
\begin{remark}\label{Rem-Short-Hand}
{Hereafter, we use the short-hand notations}:
\begin{equation}\label{g-w-z}
g := g(\tau)  , \  w :=w(\tau) , \ z := z(\tau) \ \
{\rm{and}} \ \ V_n :=V_n(\tau), \ U_n :=U_n(\tau) \ .
\end{equation}
\end{remark}
%%%%%%%%%%%%%%%%%%%%%%%%%%%%%%%%%%%%%%%%%%%%%%%%%%%%%%%%%%%%%%%%%%%%%%%%%%%%%%%%%%%%%%%%%%%%%%%%%%%%%%%%%%%%%%

The paper is organised as follows. In Section \ref{HDIS}, we give a quite explicit description of Hamiltonian dynamics of
the total non-autonomous system $\mathcal{S}+\mathcal{C}$ due to repeated interactions. Since our system is boson,
it is a unity preserving $\ast$-dynamics on the CCR von Neumann algebra generated by the Weyl operators. It is
a quasi-free $W^{*}$-dynamical system. We recall in Section \ref{Weyl-CF-E} formulae for entropy of the CCR quasi-free states.
In Section \ref{RepPertQFDyn}, we use them for the entropy production calculations. Section \ref{subsystem} is dedicated to
analysis of reduced dynamics of subsystems their correlations and convergence to equilibrium. We prove also a universality of
the short-time interaction limit {of this} dynamics for subsystem $\mathcal{S}$.

%%%%%%%%%%%%%%%%%%%%%%%%%%%%%%%%%%%%%%%%%%%%%%%%%%%%%%%%%
\section{{Hamiltonian Dynamics}} \label{HDIS}
%%%%%%%%%%%%%%%%%%%%%%%%%%%%%%%%%%%%%%%%%%%%%%%%%%%%%%%%%
We start this section  with analysis of evolution of the non-autonomous system $\mathcal{S}+\mathcal{C}$
with Hamiltonian (\ref{Ham-Model}).

A well-known way to avoid the problem of evolution of unbounded creation-annihilation operators is to construct
dynamics of the subsystem $\mathcal{S}$ on the unital CCR $C^*$-algebra $\mathscr{A}(\mathscr{F})$. Here
$\mathscr{A}(\mathscr{F})$ is generated on the Fock space $\mathscr{F}$ as the operator-norm closure of the
linear span $\mathscr{A}_w$ of the Weyl operator system:
\begin{equation}\label{S-Weyl}
 \{\widehat{w}(\alpha) = e^{i{\Phi(\alpha)}/\sqrt 2}\}_{\alpha \in \C} \ .
\end{equation}
Here $\Phi(\alpha):=\bar \alpha a + \alpha a^*$ is a self-adjoint operator with domain in $\mathscr{F}$ and
the CCR take then the Weyl form:
\begin{equation}\label{W-CCR}
 \widehat{w}(\alpha_1) \widehat{w}(\alpha_2) = e^{- {i} \, {\rm{Im}}(\bar{\alpha}_{1} \alpha_{2})/2} \
 \widehat{w}(\alpha_1 + \alpha_2) \ ,
 \quad \alpha_1,  \alpha_2 \in \C \ .
\end{equation}

Note that $\mathscr{A}(\mathscr{F})$ is a minimal $C^*$-algebra, which contains the span $\mathscr{A}_w$ of the Weyl operator
system (\ref{S-Weyl}). Algebra $\mathscr{A}(\mathscr{F})$ is contained in the $C^*$-algebra $\mathcal{L}(\mathscr{F})$ of
all bounded operators on $\mathscr{F}$.
%$\mathscr{A}(\mathscr{F}){\subset} \mathcal{L}(\mathscr{F})$.

Similarly we define the Weyl $C^*$-algebra $\mathscr{A(H)}\, {\subset} \, \mathcal{L}(\mathscr{H})$ over
$\mathscr{H}: = \mathscr{H}^{(N)}$ (\ref{H-space}). It is appropriate for description the system $\mathcal{S}+\mathcal{C}$,
see Remark \ref{Rem2}(ii). This algebra is generated by operators
\begin{equation}\label{Wz}
W(\zeta) = \bigotimes_{k=0}^N \widehat{w}(\zeta_k), \qquad
\zeta = \{ \zeta_k \}_{k=0}^N \in \C^{N + 1} \ , \ N \geq 1 .
\end{equation}
Using definitions of the boson operators  $\{b_j, b_{j}^*\}_{j=1}^{N}$ (\ref{S-bosons1}) and of the
sesqui-linear forms
\begin{equation}\label{b-forms}
      \langle \zeta, b\rangle  := \sum_{j=0}^N \bar{\zeta}_j b_j,  \qquad
      \langle b, \zeta\rangle  := \sum_{j=0}^N  \zeta_j b^*_j\ ,
\end{equation}
the Weyl operators (\ref{Wz}) can be rewritten as
\begin{equation}\label{Wz-bis}
W(\zeta) = \exp[i {\big(\langle \zeta, b\rangle + \langle b, \zeta\rangle\big)}/\sqrt 2] \ .
\end{equation}

We denote by $\mathfrak{C}_{1}(\mathscr{F})\subset \mathcal{L}(\mathscr{F})$, {the set of all} trace-class operators on
$\mathscr{F}$. A self-adjoint, non-negative operator $\rho\in \mathfrak{C}_{1}(\mathscr{F})$ with \textit{unit} trace is called
\textit{density matrix}. Note that the state $\omega_{\rho}(\cdot)$ generated by $\rho$ on the $C^{\ast}$-algebra of bounded
operators $\mathcal{L}(\mathscr{F})$:
\begin{equation}\label{state-S}
\omega_{\rho}(A) : = \Tr_{\mathscr{F}} (\rho \, A ) \ , \quad A \in \mathcal{L}(\mathscr{F}) \ ,
\end{equation}
is a \textit{normal} state \cite{BR1}.
In particular, (\ref{state-S}) defines a normal state on the Weyl algebra $\mathscr{A(F)}$.

Let $\{\rho_k\}_{k = 0}^N$ be collection of density matrices on $\mathscr{F}$. Then we can define normal
\textit{product-state}
\begin{equation}\label{state-S-C}
\omega_{\rho^{\otimes}}(\cdot):= \Tr_{\mathscr{H}} (\rho^{\otimes} \ \cdot \ ) \, , \quad
\rho^{\otimes}: = \otimes_{k=0}^{N} \rho_k \, ,
\end{equation}
on the $C^{\ast}$-algebra {$\mathscr{A}(\mathscr{H})$}, which is isometrically isomorphic to the tensor product
$\otimes_{k=0}^{N} \mathscr{A(F)}$ of identical $C^{\ast}$-algebras $\mathscr{A(F)}$. If we put
\begin{equation}\label{C-k}
C_k(\alpha) := \Tr_{\mathscr{F}}[\rho_k \, \widehat{w}(\alpha)] \ , \ \alpha \in \C \ ,
\end{equation}
then by (\ref{Wz}) and (\ref{C-k}) one obtains for $\rho^{\otimes}$ (\ref{state-S-C}) the representation:
\begin{equation}\label{Expect-W}
      \omega_{\rho^{\otimes}}(W(\zeta)):= \Tr_{\mathscr{H}}[\rho^{\otimes} \, W(\zeta)]
        = \prod_{k=0}^N C_k(\zeta_k) \ .
\end{equation}
%%%%%%%%%%%%%%%%%%%%%%%%%%%%%%%%%%%%%%%%%%%%%

Let $\varrho \in  \mathfrak{C}_{1}(\mathscr{H})$ be a (\textit{general}) density matrix on $\mathscr{H}$.
Then for the system $\mathcal{S}+\mathcal{C}$, the Hamiltonian dynamics $T_t : \varrho \mapsto \varrho(t)$ of initial density
matrix $\varrho(0) := \varrho $ is {defined as a unique} solution of the Cauchy problem for the {\textit{non-autonomous}}
Liouville equation
\begin{equation}\label{Cauchy-probl}
\partial_{t}\varrho(t) = L(t)(\varrho(t)) \ , \
\varrho(t)\big|_{t=0}  = \varrho \in \mathfrak{D} \subseteq \mathfrak{C}_{1}(\mathscr{H}) \ .
\end{equation}
Here $\mathfrak{D}$ denotes a suitable class of initial conditions. The time-dependent generator $L(t)$ with
${\rm{dom}}(L(t)) \subseteq \mathfrak{C}_{1}(\mathscr{H})$ is defined on the interval $[0, \tau N)$ by
(\ref{Ham-Model})-(\ref{H-n-commut}):
\begin{equation}\label{L-Gen-t}
L(t)(\varrho(t)):= -i \, [H_{N}(t),\varrho(t)] \ , \ t\in [0, \tau N) \ .
\end{equation}
The solution of the problem (\ref{Cauchy-probl}) is trace-norm ($\|\cdot\|_{1}$) differentiable family
{$\{ T_{t}(\varrho)\}_{t\in[0, N\tau )}$}.
%%%%%%%
By virtue of (\ref{Ham-n}) and (\ref{L-Gen-t}), equation (\ref{Cauchy-probl}) is \textit {autonomous} for each
of the interval
$[(n-1)\tau, n\tau)$ :
\begin{equation}\label{L-Gen-n}
L(t)(\cdot)= L_{n}(\cdot) = -i[H_n \, , \, \cdot \, ] \ ,  \ \
{\rm{}} \ \ t\in[(n-1)\tau, n\tau) \ ,  \  n \geqslant 1 \ ,
\end{equation}
i.e., the Liouvillian generator is piecewise constant (time-independent).

Since any $t \geqslant 0$ has the representation:
\begin{equation}\label{t}
t:=n(t)\tau + \nu(t) \ , \  n(t) := [t/\tau] \ \ {\rm{and}} \ \ \nu(t)\in[0, \tau) \ ,
\end{equation}
by Markovian independence of generators (\ref{L-Gen-n}), the $\|\cdot\|_{1}$-continuous solution of the Cauchy problem
(\ref{Cauchy-probl}) takes the iterative form:
\begin{eqnarray}\label{Solution-state-t}
&&\varrho(t)= T_{t}(\varrho):= T_{\nu(t),n}(T_{\tau,n-1}(\ldots T_{\tau,1}(\varrho)\ldots)) = \\
&& e^{-i\nu(t) H_n}e^{-i\tau H_{n-1}}\ldots e^{-i\tau H_1} \varrho \,
e^{i\tau H_1}\ldots e^{i\tau H_{n-1}}e^{i\nu(t) H_n} \ . \nonumber
\end{eqnarray}
Here $t\in [(n-1)\tau, n\tau)$ and {$n = n(t) < N$}. By the $\|\cdot\|_{1}$-continuity we obtain from
(\ref{Solution-state-t}) that
\begin{equation}\label{Solution-state-N}
\varrho(N\tau-0)=\varrho(N\tau) = T_{N\tau}(\varrho)= e^{-i\tau H_N}\ldots e^{-i\tau H_1} \varrho \,
e^{i\tau H_1}\ldots e^{i\tau H_N}.
\end{equation}
%%%%%%%%%%%%%%%%%%%%%%%%%%%%%%%%%%%%%%%%%%%%%%%%% Remark U %%%%%%%%%%%%%%%%%%%%%%%%%%%%%%%%%%%%%%%%%%%%%%%%%%%%%%%%%%%%%%%%
\begin{remark}\label{rem-U}
To bolster that (\ref{Solution-state-t}) gives a solution of the non-autonomous Cauchy problem (\ref{Cauchy-probl}) in the space
$\mathfrak{C}_{1}(\mathscr{F})$, we note that
\begin{equation}\label{U}
 U(t) = e^{-i\tau H_1} \cdots e^{-i\tau H_{n-1}}e^{-i\nu(t)H_{n(t)}} \ ,  \ \  t = n(t)\tau + \nu(t)\ ,
\end{equation}
is a one-parameter strongly continuous family of unitary operators on $\mathscr{H}$. Then $\varrho(t)= U(t) \varrho \, U^{\ast}(t)$
implies that for $t\in \mathbb{R}$  the map $T_t $ (\ref{Solution-state-t}) is trace- and positivity-preserving, such that
$t \mapsto \varrho(t)$ enjoys continuity in the weak operator topology. Since $\|\varrho(t)\|_1 =1$, the weak continuity implies
the $\| \cdot \|_1$-continuity of $t\mapsto \varrho(t)$, see e.g. \cite{Za}, Corollary 2.66. Hence, the map $T_t$ is a
trace-norm continuous $\ast$-automorphism of the set of all density matrices:
$\{\varrho \in \mathfrak{C}_{1}(\mathscr{F}): \varrho \geq 0,  \|\varrho\|_1 =1\}$.
\end{remark}
%%%%%%%%%%%%%%%%%%%%%%%%%%%%%%%%%%%%%%%%%%%%%%% Remark Dual %%%%%%%%%%%%%%%%%%%%%%%%%%%%%%%%%%%%%%%%%%%%%%%%%%%%%%%%%%%%%%%%%
\begin{remark}\label{rem-dual}
Note that equivalent and often more convenient description of density matrices evolution (\ref{Solution-state-t})
is the \textit{dual} {dynamics $T_t^{\ast} : \mathcal{L}(\mathscr{H}) \to \mathcal{L}(\mathscr{H}) $}:
\begin{equation}\label{dual}
\omega_{\, T_{t}(\varrho)}(A) = \Tr_{\mathscr{H}} (T_{t}(\varrho) \, A)
=: \Tr_{\mathscr{H}} ( \varrho \, T_{t}^{{\ast}}(A)) , \ {\rm{for}}\,
(\varrho, A) \in \mathfrak{C}_1(\mathscr{F})\times\mathcal{L}(\mathscr{H}).
\end{equation}
Since $t\mapsto T_{t}(\varrho)$ is $\|\cdot\|_{1}$-continuous and since $\mathcal{L}(\mathscr{H})$ is topologically dual of
$\mathfrak{C}_{1}(\mathscr{H})$, one gets that $t \mapsto T_{t}^{\ast}(A)$ in (\ref{dual}) is a one-parameter $\ast$-automorphism
of the unital $C^*$-algebra of bounded operators $\mathcal{L}(\mathscr{H})$. The automorphism of the $C^*$-dynamical system
$(\mathcal{L}(\mathscr{H}), T_{t}^{\ast})$ is not time-continuous for bosons, see Appendix \ref{App1}.
To ensure the continuity of $T_{t}^{\ast}$ one considers instead of the $C^*$-algebra
$\mathscr{A(H)} \subset \mathcal{L}(\mathscr{H})$, the von Neumann algebra $\mathfrak{M}(\mathscr{H})$, which is closure
of the Weyl linear span {$\mathscr{A}_w$} generated by (\ref{S-Weyl}), (\ref{Wz}) in the weak*-topology.
Since it is weaker than $C^*$-algebra topology, $\mathfrak{M}(\mathscr{H})$ is $\ast$-isomorphic to  $\mathcal{L}(\mathscr{H})$.
Then $\|\cdot\|_{1}$-continuity of $T_{t}(\varrho)$ implies continuity of the dual mapping $t \mapsto T_{t}^{\ast}(A)$ in the
weak*-topology on $\mathfrak{M}(\mathscr{H})$ and defines a $W^*$-dynamical system
$(\mathfrak{M}(\mathscr{H}), T_{t}^{\ast})$.
\end{remark}

We leave details for Appendix \ref{App1} and we address the readers to the references \cite{AJP1} and \cite{BR1}, \cite{BR2}.
%\end{remark}
%%%%%%%%%%%%%%%%%%%%%%%%%%%%%%%%%%%%%%%%%%%%%%% Remark Explicit %%%%%%%%%%%%%%%%%%%%%%%%%%%%%%%%%%%%%%%%%%%%%%%%%%%%%%%%%%%%%%%%%
\begin{remark}\label{rem-dual-explicit}
Below we show that $T_t^{\ast}$ maps $\mathscr{A}(\mathscr{H})$ into itself, and that
the action of $T_t^{\ast}$ on Weyl operators can be calculated in the explicit form.
Since $\mathscr{A}(\mathscr{H})$ is $\ast$-weakly dense in $\mathcal{L}(\mathscr{H})$, these allow to deduce properties
of evolution $\rho(t)$.
%%%%%%%%%%%%%%
%We leave details for Appendix \ref{App1} and we address the readers to the references \cite{AJP1} and \cite{BR1}, \cite{BR2}.
\end{remark}
%%%%%%%%%%%%%%%%%%%%%%%%%%%%%%%%%%%%%%%%%%%%%%%%%%%%%%%%%%%%%%%%%%%%%%%%%%%%%%%%%%%%%%%%%%%%%%%%%%%%%%%%%%%%%%%%%%%%

Using (\ref{Solution-state-N}) and dual representation (\ref{dual}), we can prove the main result of this section.
%%%%%%%%%%%%%%%%%%%%%%%%%%%%%%%%%%%%%%%%%%%%%%%%%%%%% Lemma %%%%%%%%%%%%%%%%%%%%%%%%%%%%%%%%%%%%%%%%%%%%%%%%%%%%%
\begin{lem}\label{lem21}
For $t = N\tau$, the expectation (\ref{Expect-W}) of the Weyl operator (\ref{Wz-bis}) with respect to the
evolved state has the form
\begin{equation}\label{Expect-Wt}
\omega_{\rho(N\tau)}(W(\zeta)) = \omega_{\rho}(W(U_1\ldots U_N \ \zeta))
= \prod_{k=0}^NC_k((U_1\ldots U_N \ \zeta)_k) \, .
\end{equation}
Here
\begin{equation}\label{U0}
      (U_1\ldots U_N \ \zeta)_{0} = e^{iN\tau\epsilon}\big((gz)^{N}\zeta_0
       + \sum_{j=1}^N gw(gz)^{j-1} \ \zeta_j \big) \, ,
\end{equation}
whereas
\begin{equation}\label{Uk}
(U_1\ldots U_N \ \zeta)_k = e^{iN\tau\epsilon}\big(g w (gz)^{N-k}\zeta_0 + g\bar z\zeta_k +
\sum_{j=k+1}^N g^2 w^2(gz)^{j-k-1}\zeta_j \big) \, ,
\end{equation}
for $0<k<N$,  and
\begin{equation}\label{UN}
(U_1\ldots U_N \ \zeta)_{N} = e^{iN\tau\epsilon}\big(g w \zeta_0 + g\bar z\zeta_N \big) \, ,
\end{equation}
see definitions (\ref{g-w}) and (\ref{z}).
\end{lem}
%%%%%%%%%%%%%%%%%%%%%%%%%%%%%%%%%%%%%%%%%%%%%%%%%%% Proof %%%%%%%%%%%%%%%%%%%%%%%%%%%%%%%%%%%%%%%%%%%%%%%%%%%%%%%%%%%
{\sl Proof :} Note that (\ref{Expect-W}), (\ref{Solution-state-N}) and the duality (\ref{dual}) yield
\[
        \omega_{\rho(N\tau)}(W(\zeta))= \Tr_{\mathscr{H}}[\rho \ T_{N\tau}^{{\ast}}(W(\zeta))]
        = \Tr_{\mathscr{H}}[\rho \ e^{i\tau H_1}\ldots e^{i\tau H_N}W(\zeta)
              e^{-i\tau H_N}\ldots e^{-i\tau H_1}]
\]
\begin{equation}\label{UU}
         = \Tr_{\mathscr{H}}  [\rho \ W(U_1\ldots U_N \ \zeta)]
        = \prod_{k=0}^NC_k((U_1\ldots U_N \ \zeta)_k).
\end{equation}
To {generate} the mapping $\zeta \mapsto U_1\ldots U_N \, \zeta$ in (\ref{UU}),
we use Lemma \ref{lem11} and sesquilinear forms (\ref{b-forms}) to obtain
\begin{eqnarray}\label{timeevo1}
&&e^{i\tau H_1}\ldots e^{i\tau H_N} \langle\zeta, b\rangle e^{-i\tau H_N}\ldots e^{-i\tau H_1} =
\langle\zeta,U_N^{\ast} \ldots U_1^{\ast} b\rangle  \\
&& \hspace{3cm}=\langle U_1\ldots U_N \ \zeta, b\rangle  \nonumber
\end{eqnarray}
and the similar expression for its conjugate, which we then insert into (\ref{Wz-bis}).

Moreover, by the same Lemma \ref{lem11},  we get that $U_1\ldots U_N \ \zeta = e^{iN\tau\epsilon}V_1\ldots V_N \, \zeta$,
where
\[
       (V_1\ldots V_N)_{0j} =
      \begin{cases} (V_{1})_{00}\ldots (V_N)_{00} = (gz)^{N}
                                       & \quad (j = 0) \\
      (V_{1})_{00}\ldots (V_{j-1})_{00}(V_j)_{0j}
    (V_{j+1})_{jj}\ldots(V_N)_{jj} =  (gz)^{j-1} gw
                                        & \quad (0 < j \leqslant N ) ,
     \end{cases}
\]
and for $ 0 < k  \leqslant N$:
\[
  (V_1\ldots V_N)_{kj} =
      \begin{cases} (V_1 \ldots V_{k-1})_{kk}
       (V_k)_{k0}(V_{k+1}\ldots V_N)_{00} = g w(gz)^{N-k}
                                       &  ( j = 0 ) \\
                    0  &  (0< j < k  )                   \\
       (V_1\ldots V_{k-1})_{kk} (V_k)_{kk}
            (V_{k+1}\ldots V_N)_{kk}= g\bar z &   (j = k) \\
      (V_1\ldots V_{k-1})_{kk} (V_k)_{k0}
         (V_{k+1}\ldots V_{j-1})_{00}(V_j)_{0j}
       (V_{j+1}\ldots V_N)_{jj} \hskip-5mm &           \\
         \qquad \qquad \qquad \qquad = gw(gz)^{j-k-1} gw
                                        &  (k < j \leqslant N).
     \end{cases}
\]
Collecting these formulae, one obtains explicit expressions for components  (\ref{U0}) and (\ref{Uk}) of the
vector $U_1\ldots U_N \, \zeta$.  \hfill $\square$
%%%%%%%%%%%%%%%%%%%%%%%%%%%%%%%%%%%%%%%%%%%%%%%%%%%%% Remark %%%%%%%%%%%%%%%%%%%%%%%%%%%%%%%%%%%%%%%%%%%%%%%%%%%%%%%%%%%%%%%%
\begin{remark}\label{t-m-tau}
{Note that for a fixed $N$ and for any $t = m \tau$, $1\leq m \leq N$, the arguments of Lemma \ref{lem21} give a
general formula}
\begin{eqnarray}\nonumber
&&\omega_{\rho(m\tau)}(W(\zeta)) = \omega_{\rho}(T_{ m\tau}^{{\ast}}(W(\zeta))) = \omega_{\rho}(W(U_1\ldots U_m \ \zeta))\\
&&   \hspace{2.5cm} = \prod_{k=0}^N C_k((U_1\ldots U_m \ \zeta)_k)  \label{Expect-Wm-tau} \ .
\end{eqnarray}
{Following the same line of reasoning as for (\ref{timeevo1}) one obtains explicit formulae for the components
$\{(U_1\ldots U_m \ \zeta)_k\}_{k=0}^{N}$:}
\smallskip
\begin{equation*}
\hspace{-9.5cm}(U_1\ldots U_m \ \zeta)_k =
\end{equation*}
\[
     \begin{cases} e^{im\tau\epsilon}\big((gz)^{m}\zeta_0 + \sum_{j=1}^m gw(gz)^{j-1}
           \ \zeta_j \big)   & (k = 0) \\
           e^{im\tau\epsilon}\big(g w (gz)^{m-k}\zeta_0 + g\bar z\zeta_k
         + \sum_{j=k+1}^m g^2 w^2(gz)^{j-k-1}\zeta_j \big)
                       & ( 1 \leqslant k < m )  \\
        e^{im\tau\epsilon} \big(g w \zeta_0 + g \bar z \zeta_m \big)
                      & ( k = m ) \\
         e^{im\tau\epsilon} \ \zeta_k  & ( m < k \leqslant N)
     \end{cases}
\]
{Note that for $m=N$, these formulae coincide with (\ref{U0})-(\ref{UN}), except the last line,
which is void in this case.}
\end{remark}

%%%%%%%%%%%%%%%%%%%%%%%%%%%%%%%%%%%%%%%%%%%%%%%%%%%%%%%%%%%%%%%%%%%%%%%%%%%%%%%%%%%%%%%%%%%%%%%%%%%%%%%%%%%%%%%%%%%%%%%%%%%%%
%%%%%%%%%%%%%%%%%%%%%%%%%%%%%%%%%%%%%%%%%%%%%%%%%%%%% Remark %%%%%%%%%%%%%%%%%%%%%%%%%%%%%%%%%%%%%%%%%%%%%%%%%%%%%%%%%%%%%%%%
\begin{remark}\label{QF-dyn}
Recall that unity preserving $\ast$-dynamics $t \mapsto T_{t}^*$ on the CCR von Neumann algebra
$\mathfrak{M}(\mathscr{H})$ generated by $\{W(\zeta)\}_{\zeta \in \mathbb{C}}$  (\ref{Wz-bis})
is quasi-free, if there {exist} a mapping $U_t : \zeta \mapsto U_t \zeta$ and a complex-valued function
$\Omega_t : \zeta \mapsto \Omega_t(\zeta)$,  such that
\begin{equation}\label{Q-F-map}
 T_{t}^*(W(\zeta)) = \Omega_t(\zeta) W(U_t \zeta) \ , \ \Omega_{0} = 1 \ , \ U_{0} = I \ ,
\end{equation}
see e.g., \cite{AJP1}, \cite{BR2} or \cite{Ve}.
Then by Remark \ref{t-m-tau}, the step-wise dynamics
\begin{equation*}
T_{m\tau}^*(W(\zeta)) = W(U_1\ldots U_m \ \zeta) \  ,
\qquad { m = 0, 1, \ldots , N } \ \
\end{equation*}
is quasi-free, with  $\Omega_t(\zeta) = 1$ and the matrices $\{U_j\}_{j=1}^N$  on $\C^{N + 1}$ defined by
Lemma \ref{lem11}.
\end{remark}

%%%%%%%%%%%%%%%%%%%%%%%%%%%%%%%%%%%%%%%%%%%%%%%%%%%%%%%%%%%%%%%%%%%%
%%%%%%%%%%%%%%%%%%%%%%%%%%%%%%%%%%%%%%%%%%%%%%%%%%%%%%%%%%%%%%%%%%%%
\section{{Entropy of CCR Quasi-Free States}}\label{Weyl-CF-E}
%%%%%%%%%%%%%%%%%%%%%%%%%%%%%%%%%%%%%%%%%%%%%%%%%%%%%%%%%%%%%%%%%%%%

In this section, we establish some useful formulae relating expectations of the Weyl operators (Weyl characteristic
function) and the entropy of boson quasi-free states.
We formulate them in a way that is restricted but sufficient for our purposes.
For general settings see, e.g. \cite{Fa}, \cite{AJP1}, \cite{BR2}, \cite{Ve} and references therein.
%%%%%%%%%%%%%%%%%%%%%%%%%%%%%%%%%%%%%%%%%%%%%%%%%%%%% Remark %%%%%%%%%%%%%%%%%%%%%%%%%%%%%%%%%%%%%%%%%%%%%%%%%%%%%%%%%%%%%%%%
\begin{definition}\label{QF-state}
{A state $\omega$ on the CCR $C^*$-algebra $\mathscr{A}(\mathscr{F})$ (\ref{S-Weyl}) is called {quasi-free},
if its {characteristic} function has the form
\begin{equation}\label{Char-QF1}
\omega(\widehat{w}(\alpha)) := e^{- \frac{1}{4}|\alpha|^2 - \frac{1}{2} h(\alpha)} \ , \
\end{equation}
where $h: \alpha \mapsto \widehat{h}(\alpha, \alpha)$ is a (closable) non-negative sesquilinear form on
$\mathbb{C} \times \mathbb{C}$.
A quasi-free state {$\omega$ is} gauge-invariant if
$\omega(\widehat{w}(\alpha)) = \omega(\widehat{w}(e^{i \varphi}\alpha))$ for $\varphi \in [0, 2\pi)$. }
\end{definition}

{Let $\omega_{\beta}$ denote the Gibbs state with parameter $\beta$ (inverse temperature) given by the density matrix
$\rho(\beta) =e^{-\beta a^*a}/Z(\beta)$, where  $Z(\beta)= (1 - e^{-\beta})^{-1}$.
This state is quasi-free and gauge-invariant, since
\begin{equation}\label{Char-G1}
         \omega_{\beta}(\widehat{w}(\alpha))
        = e^{- \frac{1}{4}|\alpha|^2 - \frac{1}{2} h_{\beta}(\alpha)}
\end{equation}
holds for
\begin{equation}\label{h}
h_{\beta}(\alpha)= \frac{|\alpha|^2 }{e^\beta - 1} \ , \ \alpha \in \mathbb{C} \ .
\end{equation}
}
Note that the entropy of $\omega_{\beta} $ is given by
\begin{equation}
s(\beta) := - \Tr_{\mathscr{F}}[\rho(\beta)\ln\rho(\beta)] =
\frac{\beta}{e^{\beta}-1}-\ln(1-e^{-\beta}).
\label{1entropy}
\end{equation}
and that
\begin{equation}
\omega_{\beta} (a^*a) = \frac{1}{e^{\beta}-1}  .
\label{traa}
\end{equation}

In terms of the variable
\begin{equation}\label{x}
         \frac{1+e^{-\beta}}{1-e^{-\beta}} =: x,
\end{equation}
the entropy (\ref{1entropy}) can be represented as
\begin{equation} \label{s-beta}
s(\beta)=\sigma\Big(\frac{1+e^{-\beta}}{1-e^{-\beta}}\Big)  \ ,
\end{equation}
where
\begin{equation}\label{sigma}
         \sigma(x) : = \frac{x+1}{2}\ln\frac{x+1}{2}
           - \frac{x-1}{2}\ln\frac{x-1}{2} \ .
\end{equation}
Here $ \sigma : (1, \infty) \to (0, \infty)$ and $\sigma ' (x) > 0 $.

Extension to the case of the space (\ref{H-space}) is straightforward: general \textit{gauge-invariant} quasi-free states
on the CCR $C^*$-algebra $\mathscr{A}(\mathscr{H})$ are defined by density matrices of the form \cite{Ve}:
\begin{equation}\label{ro-L}
     \rho_L = \frac{1}{Z_L} \, e^{-{\langle b, L b\rangle}} \ , \ Z_L
     = \det \big[1-e^{-L}\big]^{-1} \ .
\end{equation}
Here sesquilinear operator-valued forms $\langle b, Lb\rangle = \sum_{n,m=0}^N \ell_{nm} b^*_n b_m$
are parameterised by $ (N+1)\times(N+1) $ \textit{positive-definite} Hermitian matrix
$ L = \{\ell_{nm}\}_{0\leqslant n, m\leqslant N}$.
Note that the $\ast$-automorphism $T_\varphi$ on $\mathscr{A}(\mathscr{H})$
(the gauge transformation) :
\begin{equation}\label{g-T}
      T_\varphi : b^*_n \mapsto b^*_n e^{i \varphi}, \ b_m
       \mapsto b_m e^{- i \varphi} \qquad  (\varphi \in \mathbb{R} \ , \
n,m = 0,1,\ldots N )\ ,
\end{equation}
leaves the state (\ref{ro-L}) invariant.

Then characteristic function of the Weyl operators $W(\zeta)$ takes the form
\begin{equation}\label{W-ro-L}
\omega_{\rho_L}(W(\zeta) ) =\Tr_{\mathscr{H}}  [\rho_L W(\zeta) ]
= \exp\big[-\frac{1}{4}\langle \zeta,\zeta\rangle - \frac{1}{2}\langle \zeta,\frac{1}{e^{L}-1}\, \zeta \rangle\big] \ .
\end{equation}
Here the vector in the argument is
\[
              \zeta = \begin{pmatrix} \zeta_0 \\ \zeta_1 \\ \cdot \\
               \cdot \\ \cdot \\ \zeta_N
                      \end{pmatrix} \in \C^{N+1} \ .
\]
Note that the entropy of the state $\omega_{\rho_L}$ is given by
\begin{equation}\label{entr-L}
       S(\rho_L) = -\Tr_{\mathscr{H}} [\rho_L \ln\rho_L] =
             \tr [L(e^L -1)^{-1} - \ln(1-e^{-L})],
\end{equation}
where the trace in the third member is over $\C^{N+1}$.

If we define the matrix
\begin{equation}\label{matr-X}
        X := (1+e^{-L})(1-e^{-L})^{-1} ,
\end{equation}
then the characteristic function (\ref{W-ro-L}) takes the form:
\begin{equation}\label{weylcorrgen}
\omega_{\rho_L}(W(\zeta))
 = \exp\big[-\frac{1}{4}\langle \zeta, X\zeta \rangle\big] \, .
\end{equation}
And for the entropy (\ref{entr-L}), we obtain
 \begin{equation}\label{entropycorrgen}
        {S(\rho_L) = \tr \Big[\frac{X+1}{2}
        \ln \frac{X+1}{2} - \frac{X-1}{2} \ln \frac{X-1}{2} \Big] .}
        \end{equation}

Below we need a bit more specified set up than (\ref{matr-X})-(\ref{entropycorrgen}).
Let $ \rho(\beta, \delta ;\xi ) $ be density matrix of a quasi-free state
(\ref{ro-L}) corresponding to the operator-valued sesquilinear form
\begin{equation}
  \langle b, L(\beta, \delta; \xi)b\rangle =
     \beta\sum_{n=0}^{N}b_n^*b_n + \delta \langle b, \xi \rangle
    \langle \xi, b \rangle .
\label{ro-delta}
\end{equation}
on $\C^{N+1}\times\C^{N+1}$. Here $  \beta >0, \; \delta > - \beta, $ and the vector
\[
          \xi = \begin{pmatrix} \xi_0 \\ \xi_1 \\ \cdot \\
               \cdot \\ \cdot \\ \xi_N
                      \end{pmatrix} \in \C^{N+1} \ ,
\]
%%%%%%%%%%%%%%%%%%%%%%%%%%%%%%%%%%%%%%%%%%%%%%%%%%%% Lemma-Remark %%%%%%%%%%%%%%%%%%%%%%%%%%%%%%%%%%%%%%%%%%%%%%%%%%%%%%
\begin{lem}\label{Rem4}
The partition function of the state $ \rho(\beta, \delta ;\xi ) $ is given by
\begin{equation}\label{Z-delta-xi}
       Z(\beta, \delta ; \xi) = \Tr_{\mathscr{H}}
     [e^{-\langle b, L(\beta, \delta ; \xi) b \rangle}]=
(1-e^{-\beta})^{-N}(1-e^{-(\beta+\delta \langle \xi, \xi \rangle)})^{-1} \ ,
\end{equation}
so that
\[
   \rho(\beta, \delta ;\xi ) = \frac{1}{Z(\beta, \delta ; \xi)}\,
       \exp \big[-\langle b, L(\beta, \delta ; \xi) b \rangle \big] \ .
\]
The characteristic function and the entropy of this state are respectively:
\[
           \Tr_{\mathscr{H}}  [ \rho(\beta, \delta; \xi) W(\zeta)  ]
       = \exp\big[-\frac{1}{4} \ \frac{1+e^{-\beta}}{1-e^{-\beta}} \,
                \langle \zeta, \zeta \rangle\big]
\]
\begin{equation}
        \times \exp\big[-\frac{1}{4}
          \Big( \frac{1+e^{-\beta-\delta\langle \xi, \xi \rangle}}{1-e^{-\beta -\delta\langle \xi, \xi \rangle}}
         - \frac{1+e^{-\beta}}{1-e^{-\beta}}\Big)
           |\langle \xi, \zeta \rangle|^2/\langle \xi, \xi \rangle \big]
\label{weyl1}
\end{equation}
and
\[
        S(\rho(\beta, \delta; \xi) ) = -\Tr_{\mathscr{H}} [\rho(\beta, \delta; \xi)
            \ln \rho(\beta, \delta; \xi) ]
\]
\begin{equation}
                 = Ns(\beta) + \ s(\beta+\delta \langle \xi, \xi \rangle ) \ .
\label{entropy1}
\end{equation}

\end{lem}
%%%%%%%%%%%%%%%%%%%%%%%%%%%%%%%%%%%%%%%%%%%%%%%%%%%%%%%%%%%%%%%%%%%%%%%%%%%%%%%%%%%%%%%%%%%%%%%%%%
{\sl Proof :} {Proof of (\ref{Z-delta-xi}) follows from (\ref{ro-L}) and (\ref{ro-delta}). Indeed, since
by (\ref{ro-L}) any orthogonal transformation $\mathcal{O}$ on $\C^{N+1}$ leaves
the partition function invariant: $Z_{\mathcal{O}^{T}L\mathcal{O}} = Z_{L}$, one can calculate it with $O\xi$
(instead of $\xi$), where $\mathcal{O}\xi$ has only one non-zero component equals to the vector norm
$\langle \xi, \xi \rangle^{1/2}$. Then the right-hand side of (\ref{Z-delta-xi}) follows straightforwardly from the
calculation of the left-hand side for this choice of $\mathcal{O}\xi$.}

{Since this transformation $\mathcal{O}$ also diagonalise the matrix $L:=L(\beta, \delta; \xi)$, one uses it to simplify the
explicit calculations in (\ref{matr-X}), (\ref{weylcorrgen}) and then returns back to $\xi$ as a last step. To this aim
note that}
\begin{eqnarray}\label{X-L}
&&\omega_{\rho_L}(W(\zeta))
 = \exp\big[-\frac{1}{4}\langle O\zeta, OXO^* O\zeta \rangle\big]  = \\
&&\exp\big[-\frac{1}{4} \ \frac{1+e^{-\beta}}{1-e^{-\beta}} \,\langle O\zeta, O\zeta \rangle^{\prime}\big]
\exp\big[-\frac{1}{4}\ \frac{1+e^{-\beta-\delta\langle \xi, \xi \rangle}}{1-e^{-\beta -\delta\langle \xi, \xi \rangle}}
|(O\zeta)_0|^2 \big] \nonumber \ .
\end{eqnarray}
Here $\langle O\zeta, O\zeta \rangle^{\prime}:= \sum_{k=1}^{N}|(O\zeta)_{k}|^2$ and
we choose transformation $O$ in such a way that $(O\xi)_{j} = \delta_{0,j} \|\xi\|$. Since
\begin{equation}\label{zeta-0}
|(O\zeta)_0|^2 = \frac{1}{\langle \xi, \xi \rangle}
\langle O\zeta , O\xi \rangle \langle O\xi, O\zeta \rangle \ ,
\end{equation}
the identities (\ref{X-L}) give the proof of (\ref{weyl1}).

The same method is valid for entropy (\ref{entr-L}). Calculation of the trace in diagonal
representation for $L=L(\beta, \delta; \xi)$ gives formula (\ref{entropy1}). \hfill $\square$
%%%%%%%%%%%%%%%%%%%%%%%%%%%%%%%%%%%%%%%%%%%%%%%%%%%%%%%%%%%%%%%%%%%%%%%%%%%%%%%%%%%%%%%%%%%%%%%%%%%%%%%%%%%%%%%

Recall that the state $\omega$ on the CCR $C^*$-algebra $\mathscr{A}(\mathscr{H})$ is \textit{regular}, if the map
$s \mapsto \omega(W(s \, \zeta))$ is a continuous function of $s\in \mathbb{R}$ for any $\zeta \in \C^{N+1}$.
This property follows from the explicit expression  (\ref{weyl1}).
Since by the Araki-Segal theorem (see e.g. \cite{AJP1} or \cite{BR1}) a regular state is completely defined by its
characteristic function, (\ref{weyl1}) and (\ref{entropy1}) yield the following statement.

%%%%%%%%%%%%%%%%%%%%%%%%%%%%%%%%%%%%%%%%%%%%%%%%%% Lemma %%%%%%%%%%%%%%%%%%%%%%%%%%%%%%%%%%%%%%%%%%%%%%%%%%%%%%%%
\begin{lem}\label{lem-W-S}
The entropy $S(\rho)$ of the quasi-free state $\omega_{\rho}$ on the CCR $C^*$-algebra $\mathscr{A}(\mathscr{H})$
with characteristic function
\begin{equation}
\omega_{\rho}(W(\zeta)) =
\exp\Big[-\frac{1}{4}\big(x\langle \zeta, \zeta \rangle
+  x_0 |\langle \xi, \zeta \rangle|^2 \big) \Big]
\label{weylcorr}
\end{equation}
is uniquely determined by the parameters $(\xi,x,x_0)$,
where $\xi \in \C^{N+1}$,

\noindent$x > 1, \ x_0 > 1-x \, $
and it has the form
\begin{equation}
          S(\rho) = N\sigma(x) + \sigma(x+x_0\langle \xi, \xi \rangle ) \ ,
\label{entropycorr}
\end{equation}
where $\sigma(\cdot)$ is defined by (\ref{sigma}).
\end{lem}
%%%%%%%%%%%%%%%%%%%%%%%%%%%%%%%%%%%%%%%%%%%%%%%%%%%%%%%%%%%%%%%%%%%%%%%%%%%%%%%%%%%%%%%%%%%%%%%%%%
{\sl Proof :} The proof follows directly from definitions (\ref{x}), (\ref{s-beta}),
if one puts
\begin{equation*}
                 x_0 \langle \xi, \xi \rangle
       = \frac{1+e^{-\beta-\delta\langle \xi, \xi \rangle }}
        {1-e^{-\beta -\delta\langle \xi, \xi \rangle }} - \frac{1+e^{-\beta}}{1-e^{-\beta}} \ ,
\end{equation*}
in (\ref{weyl1}) and use (\ref{s-beta}) in (\ref{entropy1}). \hfill $\square$

%%%%%%%%%%%%%%%%%%%%%%%%%%%%%%%%%%%%%%%%%%%%%%%%%%%%%%%%
\section{Repeated Perturbations and Entropy Production} \label{RepPertQFDyn}
%%%%%%%%%%%%%%%%%%%%%%%%%%%%%%%%%%%%%%%%%%%%%%%%%%%%%%%%
We consider evolution (\ref{Solution-state-t}) of {the} system $\mathcal{S}+\mathcal{C}$, when initial density
matrix (\ref{state-S-C}) corresponds to the product of gauge-invariant Gibbs quasi-free states with parameter
$\beta_0 \geq 0$ for $\mathcal{S}$ and with parameter $\beta \geq 0$ for $\mathcal{C}$:
\begin{equation*}
    \rho = \rho_0\otimes\bigotimes_{k=1}^N  \rho_k \ ,
\end{equation*}
{where}
\begin{equation}\label{Gibbsdensity}
\rho_0= e^{-\beta_0 a^*a}/Z(\beta_0) \ , \
\rho_j= e^{-\beta a^*a}/Z(\beta) \ , \qquad
{j = 1,\ldots, N } \ .
\end{equation}

This case corresponds to $\rho_L$ in (\ref{ro-L}) with diagonal matrix
$L = $diag $ (\beta_0, \beta, \cdots , \beta)$ and to
$\rho(\beta, \delta ; \xi ) $ in representation (\ref{ro-delta}) with
$ (\beta, \delta ; \xi) = (\beta, \beta_0-\beta; e)  $, i.e.,
\begin{equation}\label{prodG}
       \rho = \rho(\beta, \beta_0-\beta; e)
       = \exp\big[-\beta_0b_0^*b_0 - \beta\sum_{j=1}^{N}b_j^*b_j\big]
           /Z(\beta, \beta_0-\beta) \ .
\end{equation}
Here
\[
           e = \begin{pmatrix} 1 \\ 0 \\ \cdot \\
               \cdot \\ \cdot \\ 0
                      \end{pmatrix} \in \C^{N+1}
\]
and
\[
         Z(\beta, \beta_0-\beta) = Z(\beta_0) Z(\beta)^N =
            \frac{1}{(1-e^{-\beta_0})(1-e^{-\beta})^N} .
\]

By a straightforward application of formulae (\ref{weyl1}), (\ref{entropy1})
and Lemma \ref{Rem4} for $\xi =e$ (i.e. $\langle \xi,\xi \rangle = 1$,
$\langle \xi,\zeta \rangle = \zeta_0 $) to the state (\ref{Gibbsdensity}) (or (\ref{prodG})),
one obtains the proof of the following statement:
%%%%%%%%%%%%%%%%%%%%%%%%%%%%%%%%%%%%%%%%%%%%%%%%%%%%%%%%% Lemma %%%%%%%%%%%%%%%%%%%%%%%%%%%%%%%%%%%%%%%%%%%%%%%%%%%%%%%
\begin{lem} \label{lem41}
The characteristic function of (\ref{Gibbsdensity}) (or (\ref{prodG})) is
 \begin{eqnarray}\label{W-L41}
&&\omega_{\rho}(W(\zeta)) = \Tr_{\mathscr{H}}[\rho \ W(\zeta)] = \\
&&\exp\Big[-\frac{|\zeta_0|^2}{4}
         \Big(\frac{1+e^{-\beta_0}}{1-e^{-\beta_0}}
       - \frac{1+e^{-\beta}}{1-e^{-\beta}}\Big)
        - \frac{\langle \zeta,\zeta \rangle}{4}
        \frac{1+e^{-\beta}}{1-e^{-\beta}}\Big] \ , \nonumber
\end{eqnarray}
and the entropy is equal to
\begin{equation}\label{S-L41}
         S(\rho) = Ns(\beta) + s(\beta_0) \,.
\end{equation}
\end{lem}
%%%%%%%%%%%%%%%%%%%%%%%%%%%%%%%%%%%%%%%%%%%%%%%%%%%%%%%%%%%%%%%%%%%%%%%%%%%%%%%%%%%%%%%%%%%%%%%%%%%%%%%%%%%%%%%%%%%%%%%%

Since by (\ref{Solution-state-t}) the density matrix $\rho(t)$ for $t=N\tau$ is
       \begin{equation}\label{ro-N-tau}
       \rho(N\tau) = e^{-i\tau H_N}
    \dots e^{-i\tau H_1} \rho \,
       e^{i\tau H_1} \dots e^{i\tau H_N} ,
       \end{equation}
we obtain for evolution of the characteristic function and the entropy of the total system $\mathcal{S}+\mathcal{C}$:
%%%%%%%%%%%%%%%%%%%%%%%%%%%%%%%%%%%%%%%%%%%%%%%%%%%%%%%% Lemma %%%%%%%%%%%%%%%%%%%%%%%%%%%%%%%%%%%%%%%%%%%%%%%%%%%%%%%%%
\begin{lem}\label{lem42}
Characteristic function of the state with density matrix (\ref{ro-N-tau}) is equal to
\begin{eqnarray}\label{W-L42}
&&\omega_{\rho(N\tau)}(W(\zeta)) = \\
&&\exp\Big[-\frac{|(U_1 \ldots U_N\zeta)_0|^2}{4} \Big(\frac{1+e^{-\beta_0}}{1-e^{-\beta_0}}
- \frac{1+e^{-\beta}}{1-e^{-\beta}}\Big) - \frac{\langle \zeta, \zeta \rangle}{4}
\frac{1+e^{-\beta}}{1-e^{-\beta}}\Big] \ , \nonumber
\end{eqnarray}
whereas the total entropy rests invariant:
\[
        S(\rho(N\tau))  = S(\rho) = Ns(\beta) + s(\beta_0) \ .
\]
Here the mapping $U_1 \ldots U_N: \C^{N+1} \rightarrow \C^{N+1}$ is given by (\ref{U0}) and (\ref{Uk}).
\end{lem}
%%%%%%%%%%%%%%%%%%%%%%%%%%%%%%%%%%%%%%%%%%%%%%%%%%%%%%%%%%%%%%%%%%%%%%%%%%%%%%%%%%%%%%%%%%%%%%%%%%%%%%%%%%%%%%%%%%%%%%%%%
{\sl Proof }: From (\ref{Expect-Wt}), one gets
$\omega_{\rho(N\tau)}(W(\zeta)) = \omega_{\rho}(W(U_1\ldots U_N \ \zeta))$.
Since the mappings $U_j: \C^{N+1} \rightarrow \C^{N+1}$, $j=1, \ldots, N$ are
unitary (Lemma \ref{lem21}), (\ref{W-L41}) yields (\ref{W-L42}).
Finally, one gets that the mapping (\ref{ro-N-tau}) leaves the total entropy (\ref{S-L41})
invariant by definition (\ref{1entropy}) . \hfill$\square$
%%%%%%%%%%%%%%%%%%%%%%%%%%%%%%%%%%%%%%%%%%%%%%%%%%%%%%%%%%%%%%%%%%%%%%%

\medskip
%%%%%%%%%%%%%%%%%%%%%%%%%%%%%%%%%%%%%%%%%%%%%%%%%%%%%%%%%%%%%%%%%%%%%%%
Let $\omega$ and $\omega_0$ be two normal states on the Weyl CCR algebra $\mathscr{A(H)}$ with density matrices $\varrho$
and $\varrho_0$. Following Araki \cite{Ar1}, we introduce the \textit{relative entropy} of the state $\omega$ with respect to
$\omega_0$:
\begin{equation}\label{rel-S}
        {\rm{Ent}}(\varrho|\varrho_0) :=
        \Tr_{\mathscr{H}} [\varrho (\ln \varrho - \ln \varrho_0)]  \geq 0 \ ,
\end{equation}
see also \cite{AJP3}.
%%%%%%%%%%%%%%%%%%%%%%%%%%%%%%%%%%%%%%%%%%%%%%%%%%%%%%%%%% Lemma %%%%%%%%%%%%%%%%%%%%%%%%%%%%%%%%%%%%%%%%%%%%%%%%%%%%%%%%
\begin{lem} \label{lem43}
The relative entropy of $\omega_{\rho(N\tau)}$ with respect to $\omega_{\rho}$ is
\begin{equation}\label{rel-S-N}
      {\rm{Ent}}(\rho(N\tau)|\rho) \;
        = \; \frac{(\beta_0-\beta)(e^{\beta_0}-e^{\beta})}
          {(e^{\beta_0}-1)(e^{\beta}-1)}(1-|z|^{2N})\ ,
\end{equation}
where $z:=z(\tau)$ is defined by (\ref{z}) and (\ref{g-w-z}).
\end{lem}
%%%%%%%%%%%%%%%%%%%%%%%%%%%%%%%%%%%%%%%%%%%%%%%%%%%%%%%%%%%%%%%%%%%%%%%%%%%%%%%%%%%%%%%%%%%%%%%%%%%%%%%%%%%%%%%%%%%%%%%%%
{\sl Proof }: The trace cyclicity yields
\begin{equation}\label{S-rel-N}
   {\rm{Ent}}(\rho(N\tau)|\rho) =  \Tr_{\mathscr{H}}[\rho(N\tau)
     (\ln \rho(N\tau) - \ln\rho)]
\end{equation}
\[
         =  \Tr_{\mathscr{H}}[\rho
     (\ln \rho -   e^{i\tau H_1} \ldots e^{i\tau H_N}
  \ln\rho e^{-i\tau H_N}\ldots e^{-i\tau H_1} )]
\]
\[
     =  \frac{\beta -\beta_0}{Z(\beta, \beta_0 -\beta)}
      \Tr_{\mathscr{H}}\big[e^{-\beta_0b_0^*b_0 -\beta\sum_{j=1}^Nb^*_jb_j}
      \big(b^*_0b_0 - e^{i\tau H_1} \ldots e^{i\tau H_N}
  b_0^* b_0 e^{-i\tau H_N}\ldots e^{-i\tau H_1} )]\ .
\]
Note that one gets $b_0^* b_0 = \langle b, e\rangle \langle e, b\rangle $ by (\ref{b-forms}).
Hence, (\ref{timeevo1})
implies
\begin{equation}\label{b-e-e-b}
e^{i\tau H_1} \cdots e^{i\tau H_N}
  b_0^* b_0 e^{-i\tau H_N}\cdots e^{-i\tau H_1} = \sum_{k=0}^N {(U_1 \ldots U_N \ e)}_k b^*_k \ \
\sum_{k'=0}^N \overline{(U_1 \ldots U_N \ e)}_{k'} \, b_{k'}
\end{equation}

%%%%%%%%%%%%%%%%%%%%%%%%%%%%%%%%%%%%%%%%%%%%%%%%%%%%%%%%%%%%%%%%%%%%%%%%%%%%%%%%%%%%%%%%%%%%%%%%%%%%%%%%%
Note {also} that for any $k= 0,1, \ldots ,N$, one has $[ b_k^* b_k , \rho ] =0$, which implies the selection rule:
\begin{equation}\label{sel-rule}
\frac{1}{Z(\beta, \beta_0 -\beta)}
\Tr_{\mathscr{H}}\big[e^{-\beta_0b_0^*b_0 -\beta\sum_{j=1}^Nb^*_jb_j}
      \ b^*_k b_{k^\prime} \ ]  = 0 \ \ {\rm{for}} \ \ k \neq k^\prime \ .
\end{equation}
This rule implies that after injection of (\ref{b-e-e-b}) into (\ref{S-rel-N}) only diagonal terms with $k = k^\prime$
will survive in the expectation:
\begin{eqnarray*}
&& {\rm{Ent}}(\rho(N\tau)|\rho) = \\
&&  \frac{\beta -\beta_0}{Z(\beta, \beta_0 -\beta)}
\Tr_{\mathscr{H}}\big[e^{-\beta_0 b_0^*b_0 -\beta \sum_{j=1}^Nb^*_jb_j}
\big(b^*_0b_0 - \sum_{k=0}^N|(U_1 \ldots U_Ne)_k|^2 b^*_kb_k \big)\big]
\end{eqnarray*}
Finally, by Lemma \ref{lem21}, (\ref{U0}), (\ref{Uk}) and by (\ref{1entropy}), (\ref{traa}), we obtain
\begin{eqnarray*}
&& {\rm{Ent}}(\rho(N\tau)|\rho) = \\
&&  \frac{\beta -\beta_0}{Z(\beta, \beta_0 -\beta)} \Tr_{\mathscr{H}}\big[e^{-\beta_0b_0^*b_0 -\beta\sum_{j=1}^Nb^*_jb_j}
\big((1-|z|^{2N})b^*_0b_0 -\sum_{k=1}^N|w|^2|z|^{2N-2k} b^*_kb_k)\big)\big]  \nonumber \\
&&=  \frac{(\beta_0-\beta)(e^{\beta_0}-e^{\beta})} {(e^{\beta_0}-1)(e^{\beta}-1)}(1-|z|^{2N}) \ , \nonumber
\end{eqnarray*}
that proves  (\ref{rel-S-N}).                                      \hfill $\square$
%%%%%%%%%%%%%%%%%%%%%%%%%%%%%%%%%%%%%%%%%%%%%%%%%%%%%%%%%%%%%%%%%%%%%%%%%%%%%%%%%%%%%%%%%%%%%%%%%%%%%%%%%%%
\begin{remark}\label{Rem5}
The relative entropy defined by (\ref{rel-S}) is non-negative.
In contrast to \textit{invariant} total entropy (Lemma \ref{lem42}),
the relative entropy (\ref{rel-S-N}) is increasing monotonically as $N \to \infty$
for $|z|<1$ (see Lemma \ref{lem11}, Remark \ref{Rem3}). It converges to the limit:
\begin{equation}\label{rel-S-lim}
\lim_{N\to\infty} {\rm{Ent}}(\rho(N\tau)|\rho) \; =  \;
(\beta-\beta_0)\left[ \frac{1}{e^{\beta_0}-1} - \frac{1}{e^{\beta}-1} \right] \geq 0\ ,
\end{equation}
which is positive for $\beta_0 \neq \beta$. The limit (\ref{rel-S-lim}) gives asymptotic amount of the entropy
production, when one starts with the initial state corresponding to (\ref{Gibbsdensity}) and then
consider $N\tau \rightarrow \infty $, see \cite{AJP3}.
\end{remark}

%%%%%%%%%%%%%%%%%%%%%%%%%%%%%%%%%%%%%%%%%%%%%%%%%%%%%%%%%%%%%%%
%%%%%%%%%%%%%%%%%%%%%%%%%%%%%%%%%%%%%%%%%%%%%%%%%%%%%%%%%%%%%%%
\section{Evolution of Subsystems }\label{subsystem}
%%%%%%%%%%%%%%%%%%%%%%%%%%%%%%%%%%%%%%%%%%%%%%%%%%%%%%%%%%%%%%%
%%%%%%%%%%%%%%%%%%%%%%%%%%%%%%%%%%%%%%%%%%%%%%%%%%%%%%%%%%%%%%%
\subsection{Convergence to Equilibrium}\label{Conv-to-Eq}
%%%%%%%%%%%%%%%%%%%%%%%%%%%%%%%%%%%%%%%%%%%%%%%%%%%%%%%%%%%%%%%
\textit{Subsystem} $\mathcal{S}$. We start with the simplest subsystem $\mathcal{S}$.
Let the {initial} state of the total system $\mathcal{S} + \mathcal{C}$ in (\ref{H-space})
be a tensor-product of the corresponding density matrices $\rho = \rho_S \otimes \rho_C$ (see Hypothesis 1).
Then for $t\geq 0$ the state $\omega_{\mathcal{S}}^{t}(\cdot)$ of the subsystem $\mathcal{S}$ is {given} on the Weyl
$C^*$-algebra $\mathscr{A}(\mathscr{H}_{0})$ by
\begin{equation}\label{st-S}
   {\omega_{\mathcal{S}}^{t}(\cdot) := \omega_{\rho(t)}(\cdot \otimes \mathbb{1})} \ .
\end{equation}
{For $\zeta = (\alpha, 0, \ldots , 0) \in \C^{N + 1} $, let us consider the Weyl operator
$W(\zeta)= \widehat{w}(\alpha)\otimes \mathbb{1}\otimes \ldots \otimes \mathbb{1}$ (\ref{Wz}).}
By virtue of (\ref{Expect-W}), (\ref{Expect-Wm-tau}) and (\ref{st-S}){,} we obtain
for $t=m\tau$ ( $1\leq m \leq N$ ):
\begin{equation}\label{st-S-t}
\omega_{\mathcal{S}}^{m\tau}( \widehat{w}(\alpha)) = \omega_{\rho(m\tau)}(W(\zeta)) = \omega_{\rho}(W(U_1\ldots U_m \ \zeta)) \ .
\end{equation}
Then for components $\{(U_1\ldots U_m \ \zeta)_k\}_{k=0}^{N}$
of the vector $U_1\ldots U_m \, \zeta$ in (\ref{st-S-t}), one obtains the expression:
{\begin{equation}\label{0Uk0}
(U_1\ldots U_m \ \zeta)_{k} =  \begin{cases}
          e^{im\tau\epsilon} (gz)^{m}\alpha  \ & \ (k=0)  \\
          e^{im\tau\epsilon} g w (gz)^{m-k}\alpha  \ & \ (1\leqslant k < m) \\
         e^{im\tau\epsilon} g w \alpha \ & \  ( k = m )    \\
           0 \ & \ (m < k \leqslant N ) \ ,
\end{cases}
\end{equation}}
which follows from  Remark \ref{t-m-tau}.

If the initial density matrices: $\rho = \rho_S \otimes \rho_C$ corresponds to the product of Gibbs quasi-free states
for {different temperatures} as in (\ref{Gibbsdensity}), then (\ref{st-S-t}) and
Lemma\ref{lem41} yield
\begin{equation}\label{Expect-W-alpha}
\omega_{\mathcal{S}}^{m\tau}( \widehat{w}(\alpha)) =
\exp\Big[- \frac{|\alpha|^2}{4}\frac{1+e^{-\beta}}{1-e^{-\beta}}
- \frac{|z^m \alpha|^2}{4}\left(\frac{1+e^{-\beta_0}}{1-e^{-\beta_0}}
- \frac{1+e^{-\beta}}{1-e^{-\beta}}\right)\Big]
\end{equation}

Note that for any moment $t= m\tau$ the state $\omega_{\mathcal{S}}^{m\tau}(\cdot)$ is
a quasi-free Gibbs equilibrium state with parameter $\beta^*(m\tau)$ which satisfies the equation
\begin{equation}\label{*}
\frac{1+e^{-\beta^*(m\tau)}}{1-e^{-\beta^*(m\tau)}} =
          |z|^{2m}\frac{1+e^{-\beta_0}}{1-e^{-\beta_0}} +
        ( 1- |z|^{2m})\frac{1+e^{-\beta}}{1-e^{-\beta}} \ .
\end{equation}
This equation yields that either $\beta \leq \beta^*(m\tau)\leq \beta_0$, or $\beta_0 \leq \beta^*(m\tau)\leq \beta$.

The Hypothesis 5 implies that for $m \rightarrow\infty$ ($N \rightarrow\infty$) the Weyl
characteristic function (\ref{Expect-W-alpha}) has the limit
\begin{equation}\label{st-S-t=inf}
\lim_{m \rightarrow\infty}\omega_{\mathcal{S}}^{m\tau}( \widehat{w}(\alpha)) =
\exp\Big[- \frac{|\alpha|^2}{4} \frac{1+e^{-\beta}}{1-e^{-\beta}}\Big] \ .
\end{equation}
Therefore, in the limit $t \rightarrow \infty$ the subsystem $\mathcal{S}$ evolves from the Gibbs equilibrium state
with parameter $\beta_0$ to another equilibrium state with parameter $\beta$ imposed by the chain $\mathcal{C}$.

%%%%%%%%%%%%%%%%%%%%%%%%%%%%%%%%%%%%%%%%%%%%%%%%%%%%% Remark %%%%%%%%%%%%%%%%%%%%%%%%%%%%%%%%%%%%%%%%%%%%%%%%%%%%%%%%
%%%%%%%%%%%%%%%%%%%%%%%%%%%%%%%%%%%%%%%%%%%%%%%%%%%%%%%%%%

\noindent
\textit{Subsystem} $\mathcal{S}_1$. The initial state
{ $\omega_{\mathcal{S}_1}^{0}(\cdot) = \omega_{\mathcal{S}_1}^{t}(\cdot)|_{t=0}$ }
of this subsystem again corresponds to a one-point reduced density matrix or the partial trace on the CCR Weyl algebra
$\mathscr{A}(\mathscr{H}_{1})$:
\begin{equation}\label{st-S1}
\omega_{\mathcal{S}_1}^{0}(\widehat{w}(\alpha))=
\omega_{\rho}(\mathbb{1}\otimes \widehat{w}(\alpha) \otimes \bigotimes_{k=2}^{N}\mathbb{1}) =
\exp\Big[- \frac{|\alpha|^2}{4}\frac{1+e^{-\beta}}{1-e^{-\beta}}\Big]  \ .
\end{equation}
Now we choose vector $\zeta^{(1)} := (0, \alpha, 0, \ldots , 0) \in \C^{N + 1} $. Then
\begin{equation}\label{st-S1-t}
\omega_{\mathcal{S}_1}^{m\tau}( \widehat{w}(\alpha)) =
\omega_{\rho(m\tau)}(W(\zeta^{(1)}))
= \omega_{\rho_S \otimes \rho_C}(W(U_1\ldots U_m \ \zeta^{(1)}))
\end{equation}
for $1 < m \leq N$.  {By} Remark \ref{t-m-tau}{,} the components $\{(U_1\ldots U_m \ \zeta^{(1)})_k\}_{k=0}^{N}$ are:
{\begin{equation}\label{1Uk0}
(U_1\ldots U_m \ \zeta)_{k} =  \begin{cases}
          e^{im\tau\epsilon} gw \, \alpha  \ & \ (k=0)  \\
          e^{im\tau\epsilon} \, \delta_{k,1} \, g\overline{z} \, \alpha & \ (1 \leqslant k <m )\\
       0 \ & \  ( k = m )    \\
           0 \ & \ (m < k \leqslant N ) .
\end{cases}
\end{equation}}
{Then, we have}
\begin{equation}\label{Expect-W-alpha1}
\omega_{\mathcal{S}_{1}}^{m\tau}( \widehat{w}(\alpha)) =
\exp\Big[- \frac{|\alpha|^2}{4}\frac{1+e^{-\beta}}{1-e^{-\beta}}
- \frac{|w \alpha|^2}{4}\left(\frac{1+e^{-\beta_0}}{1-e^{-\beta_0}} - \frac{1+e^{-\beta}}{1-e^{-\beta}}\right)\Big]
\end{equation}
for any $1 < m \leq N$. Therefore, the initial state (\ref{st-S1}) changes to (\ref{Expect-W-alpha1}) after the first
act of interaction {on the interval $[0,\tau)$} and there is no further evolution of this state for $t>\tau$.

Note that (\ref{Expect-W-alpha1}) is characteristic function of a quasi-free Gibbs equilibrium state with parameter $\beta^*$,
which satisfies the equation
\[
     \frac{1+e^{-\beta^*}}{1-e^{-\beta^*}} =
          |w|^{2}\frac{1+e^{-\beta_0}}{1-e^{-\beta_0}} +
        ( 1- |w|^{2})\frac{1+e^{-\beta}}{1-e^{-\beta}} \ .
\]
Again, this equation implies that either $\beta \leq \beta^* \leq \beta_0$, or $\beta_0 \leq \beta^* \leq \beta$.

Evolution of subsystem $\mathcal{S}_1$ has a transparent interpretation: after the one act of interaction during the
time $t \in [0,\tau)$, it relaxed to an \textit{intermediate} equilibrium with the subsystem $\mathcal{S}$. This results in a
shift of initial parameter $\beta$ to $\beta^*$, which  rests unchangeable since there is no perturbation of $\mathcal{S}_1$
for $t> \tau$.
%%%%%%%%%%%%%%%%%%%%%%%%%%%%%%%%%%%%%%%%%%%%%%%%%%%%%%%%%%%%%%%%%%%%%%%%%%%%%%%%%%%%%%%%%%%%%%%%%%%%%%%%%%%%%%%%%%%
\smallskip

\noindent
\textit{Subsystem} $\mathcal{S}_m$. For $1 < m \leq N$ the initial state
{$ \omega_{\mathcal{S}_m}^{0}(\cdot) =  \omega_{\mathcal{S}_m}^{t}(\cdot)|_{t=0}$ }
of this subsystem is defined by the partial trace on the CCR Weyl algebra $\mathscr{A}(\mathscr{H}_{m})$:
\begin{equation}\label{st-Sm}
\omega_{\mathcal{S}_m}^{0}(\widehat{w}(\alpha))=
\omega_{\rho}(\bigotimes_{k=0}^{m-1}\mathbb{1}\otimes\widehat{w}(\alpha)\otimes
\bigotimes_{k=m+1}^{N}\mathbb{1})=\exp\Big[- \frac{|\alpha|^2}{4}\frac{1+e^{-\beta}}{1-e^{-\beta}}\Big]\ .
\end{equation}
Now we choose vector $\zeta^{(m)} := (0, \ldots , 0, \alpha, 0, \ldots , 0) \in \C^{N + 1} $, where $\alpha$
occupies the $m+1$ position. Consequently
\begin{equation}\label{st-Sm-t}
\omega_{\mathcal{S}_m}^{m\tau}( \widehat{w}(\alpha)) =
\omega_{\rho(m\tau)}(W(\zeta^{(m)})) = \omega_{\rho_S \otimes \rho_C}(W(U_1\ldots U_m \ \zeta^{(m)})) \ .
\end{equation}
The components $\{(U_1\ldots U_m \ \zeta^{(m)})_k\}_{k=0}^{N}$ are:
\begin{equation}\label{mUk0}
(U_1\ldots U_m \ \zeta^{{(m)}})_{k} =  \begin{cases}
          e^{im\tau\epsilon} \, gw (g{z})^{m-1}\, \alpha  \ & \ (k=0)  \\
          e^{im\tau\epsilon} \, g^2w^2  \, (g{z})^{m-k-1} \, \alpha & \ (1 \leqslant k <m )\\
       e^{im\tau\epsilon} g\overline{z}\, \alpha \ , \ & \  ( k = m )    \\
           0 \ & \ (m < k \leqslant N ) .
\end{cases}
\end{equation}
which again follows from {Remark \ref{t-m-tau}}. {Then, we obtain} for {evolution of the state} of the subsystem
$\mathcal{S}_m$:
\begin{eqnarray}\label{Expect-W-alpha-m}
&&\omega_{\mathcal{S}_{m}}^{m\tau}( \widehat{w}(\alpha)) = \\
&&\exp\Big[- \frac{|\alpha|^2}{4}\frac{1+e^{-\beta}}{1-e^{-\beta}}
- \frac{|w \alpha|^2}{4} |z|^{2(m-1)} \left(\frac{1+e^{-\beta_0}}{1-e^{-\beta_0}} -
\frac{1+e^{-\beta}}{1-e^{-\beta}}\right)\Big] \ . \nonumber
\end{eqnarray}

Note that interaction {for $ t\in [(m-1)\tau, m\tau)$} push out the subsystem $\mathcal{S}_m$ from the Gibbs equilibrium
state (\ref{st-Sm}), but {its effect attenuates for large $m$}:
\begin{equation}\label{S-m-relax}
\lim_{m\rightarrow\infty} \omega_{\mathcal{S}_{m}}^{m\tau}( \widehat{w}(\alpha)) =
\exp\Big[- \frac{|\alpha|^2}{4}\frac{1+e^{-\beta}}{1-e^{-\beta}}\Big]\ .
\end{equation}
Again, this is evolution of a quasi-free Gibbs equilibrium state with time-dependent inverse temperature parameter
$\beta^{**}(m\tau)$, which satisfies the equation
\begin{equation}\label{**}
     \frac{1+e^{-\beta^{**}(m\tau)}}{1-e^{-\beta^{**}(m\tau)}} =
          |w|^{2}|z|^{2(m-1)}\frac{1+e^{-\beta_0}}{1-e^{-\beta_0}} +
        ( 1- |w|^{2}|z|^{2(m-1)})\frac{1+e^{-\beta}}{1-e^{-\beta}} \ .
\end{equation}
As above, the value of the parameter $\beta^{**}(m\tau)$ is always between $\beta_0$ and $\beta$.

To interpret the evolution of $\mathcal{S}_m$ and the coincidence between (\ref{S-m-relax}) and (\ref{st-S-t=inf}) note that
the state of the subsystem $\mathcal{S}$ relaxes to that of initial state of the chain $\mathcal{C}$, see (\ref{st-S-t=inf}).
Therefore, after interaction of the subsystem $\mathcal{S}_m$, i.e. at the moment $t=m\tau$, its parameter $\beta^{**}(m\tau)$
has a value between $\beta$ and $\beta^{*}((m-1)\tau)$ since (\ref{*}) and (\ref{**}) yield
\begin{equation*}
\frac{1+e^{-\beta^{**}(m\tau)}}{1-e^{-\beta^{**}(m\tau)}} =
          |w|^{2}\frac{1+e^{-\beta^{*}((m-1)\tau)}}{1-e^{-\beta^{*}((m-1)\tau)}} +
        ( 1- |w|^{2})\frac{1+e^{-\beta}}{1-e^{-\beta}} \ .
\end{equation*}
As in the case $m=1$, one may convince that there is no further evolution:
$\omega_{\mathcal{S}_m}^{n\tau} = \omega_{\mathcal{S}_m}^{m\tau}$ for $ n \geqslant m$.

%%%%%%%%%%%%%%%%%%%%%%%%%%%%%%%%%%%%%%%%%%%%%%%%%%%%%%%%%%%%%%%%%%%%%%%%%%%%%%%%%%%%%%%%%%%%%%%%%%%%%%%%%%%%%%%%%%%
\smallskip

Next, we consider the composed subsystems $\mathcal{S} + \mathcal{S}_m$ and $\mathcal{S}_{m-n} + \mathcal{S}_{m}$.
Our aim is to study the eventual \textit{correlations} imposed by repeated perturbations due to $\mathcal{S}$.
%%%%%%%%%%%%%%%%%%%%%%%%%%%%%%%%%%%%%%%%%%%%%%%%%%%%%%%%%%%%%%%%%%%%%%%%%%%%%%%%%%%%%%%%%%%%%%%%%%%%%%%%%%%%%%%%%%%
\smallskip

\noindent
\textit{Subsystem} $\mathcal{S} + \mathcal{S}_m$. For $1 < m \leqslant N$ the initial state
$\omega_{\mathcal{S} + \mathcal{S}_m}^{0}(\cdot) =
\omega_{\mathcal{S} + \mathcal{S}_m}^{t}(\cdot)|_{t=0}$
of this \textit{composed} subsystem is defined by the partial trace on the Weyl $C^*$-algebra
$\mathscr{A}(\mathscr{H}_{0}\otimes \mathscr{H}_{m}) \approx \mathscr{A}(\mathscr{H}_{0})\otimes \mathscr{A}(\mathscr{H}_{m})$ by:
\begin{eqnarray}\nonumber
&&\omega_{\mathcal{S} + \mathcal{S}_m}^{0}(\widehat{w}(\alpha_0)\otimes\widehat{w}(\alpha_1)):=
\omega_{\rho}(\widehat{w}(\alpha_0)\otimes\bigotimes_{k=1}^{m-1}\mathbb{1}\otimes\widehat{w}(\alpha_1)
\otimes \bigotimes_{k=m+1}^{N}\mathbb{1})\\
&& = \exp\Big[- \frac{|\alpha_0|^2}{4}\frac{1+e^{-\beta_0}}{1-e^{-\beta_0}}\Big]
\exp\Big[- \frac{|\alpha_1|^2}{4}\frac{1+e^{-\beta}}{1-e^{-\beta}}\Big]\ . \label{st-S0-Sm}
\end{eqnarray}
This is the characteristic function of the product state corresponding to two isolated systems with different temperatures.
If one defines vector $\zeta^{(0,m)} := {(\alpha_0, 0, \ldots , 0, \alpha_1, 0, \ldots , 0)} \in \C^{N + 1} $,
where {$\alpha_1$} occupies the $m+1$ position, then
\begin{equation}\label{st-S0-Sm-t}
\omega_{\mathcal{S} +\mathcal{S}_m}^{m\tau}(\widehat{w}(\alpha_0)\otimes\widehat{w}(\alpha_1)) =
\omega_{\rho(m\tau)}(W(\zeta^{(0,m)})) = \omega_{\rho_S \otimes \rho_C}(W(U_1\ldots U_m \ \zeta^{(0,m)})) \ .
\end{equation}
The components $\{(U_1\ldots U_m \ \zeta^{(0,m)})_k\}_{k=0}^{N}$ are deduced from Remark \ref{t-m-tau}:
\begin{equation}\label{0mUk0}
(U_1\ldots U_m \ \zeta^{(0,m)})_k =
\begin{cases}  e^{im\tau\epsilon} \ (g{z})^{m-1} \ [g{z}\, \alpha_0 +  g w \, \alpha_1 ] , \ & (k=0) \\
e^{im\tau\epsilon} (g{z})^{m-k-1} g^2 [  w z \, \alpha_0 + w^2 \, \alpha_1 ], \ & (1\leqslant k < m) \\
e^{im\tau\epsilon} \ [g w \, \alpha_0 + g \overline{z}\, \alpha_1 ] , \  &  (k = m) \\
 0 \ & (m<k\leqslant N) .
\end{cases}
\end{equation}
Together with (\ref{Expect-W}), one gets
\begin{eqnarray}\label{Expect-W-alpha-0m}
&&\omega_{\mathcal{S} +\mathcal{S}_m}^{m\tau}(\widehat{w}(\alpha_0)\otimes\widehat{w}(\alpha_1))  \\
&&= \exp\Big[- \frac{1}{4} |z \alpha_0 + w \alpha_1|^2 |z|^{2(m-1)} \frac{1+e^{-\beta_0}}{1-e^{-\beta_0}}\Big]   \nonumber \\
&&\times \exp\Big[- \frac{1}{4} |z \alpha_0 + w \alpha_1|^2 (1 - |z|^{2(m-1)}) \frac{1+e^{-\beta}}{1-e^{-\beta}}\Big]
 \exp\Big[- \frac{1}{4} |w \alpha_0 + \overline{z} \alpha_1|^2  \frac{1+e^{-\beta}}{1-e^{-\beta}}\Big]
\nonumber \\
&&\longrightarrow \exp\Big[- \frac{1}{4} (|\alpha_0|^2 + | \alpha_1|^2)  \frac{1+e^{-\beta}}{1-e^{-\beta}}\Big] \nonumber
\end{eqnarray}
for $ m \to \infty$.

Hence, in this limit the composed subsystem $\mathcal{S} + \mathcal{S}_m$ evolves from the product of two quasi-free
equilibrium states (\ref{st-S0-Sm}) with different parameters $\beta_0$ and $\beta$ to the product of quasi-free
equilibrium states for the same parameter $\beta$ imposed by repeated interaction with the chain $\mathcal{C}$, when
$m \rightarrow \infty$. Interpretation of this is similar to the case \textit{Subsystem} $\mathcal{S}_m$.
%%%%%%%%%%%%%%%%%%%%%%%%%%%%%%%%%%%%%%%%%%%%%%%%%%%%%%%%%%%%%%%%%%%%%%%%%%%%%%%%%%%%%%%%%%%%%%%%%%%%%%%%%%%%%%%%%%%
\smallskip

\noindent
\textit{Subsystem} $\mathcal{S}_{m-n} + \mathcal{S}_{m}$. We suppose that $1 < (m-n) < m \leqslant N$. Then the initial state
$\omega_{\mathcal{S}_{m-n} + \mathcal{S}_{m}}^{t}(\cdot)|_{t=0}$ of this \textit{composed} subsystem is the
partial trace over the Weyl $C^*$-algebra
$\mathscr{A}(\mathscr{H}_{m-n}\otimes \mathscr{H}_{m})\approx \mathscr{A}(\mathscr{H}_{m-n})\otimes \mathscr{A}(\mathscr{H}_{m})$:
\begin{eqnarray}\label{st-Smn-Sm}
&&\omega_{\mathcal{S}_{m-n} + \mathcal{S}_{m}}^{0}(\widehat{w}(\alpha_{1})\otimes\widehat{w}(\alpha_2)):= \\
&&\omega_{\rho}(\bigotimes_{k=0}^{m-n-1}\mathbb{1}\otimes\widehat{w}(\alpha_{1})\otimes
\bigotimes_{k=m-n+1}^{m-1}\mathbb{1} \otimes \widehat{w}(\alpha_{2})\otimes \bigotimes_{k=m+1}^{N}\mathbb{1}) = \nonumber \\
&& = \exp\Big[- \frac{|\alpha_{1}|^2}{4}\frac{1+e^{-\beta}}{1-e^{-\beta}}\Big]
\exp\Big[- \frac{|\alpha_{2}|^2}{4}\frac{1+e^{-\beta}}{1-e^{-\beta}}\Big]\ . \nonumber
\end{eqnarray}
This is the characteristic function of the product state corresponding to two isolated systems with the same temperatures.

We define vector $\zeta^{(m-n, m)} := {(0, 0, \ldots , 0, \alpha_{1}, 0, \ldots , 0, \alpha_{2}, 0, \ldots, 0)} \in \C^{N + 1} $,
where $\alpha_{1}$ occupies the $m-n+1$ position,  and $\alpha_{2}$ occupies the $m+1$ position, then
\begin{eqnarray}\label{st-Smn-Sm-t}
&&\omega_{\mathcal{S}_{m-n} + \mathcal{S}_{m}}^{m\tau}(\widehat{w}(\alpha_{1})\otimes\widehat{w}(\alpha_{2})) = \\
&&\omega_{\rho(m\tau)}(W(\zeta^{({m-n},m)})) =
\omega_{\rho_S \otimes \rho_C}(W(U_1\ldots U_m \ \zeta^{({m-n},m)})) \ . \nonumber
\end{eqnarray}
Again, with help of Remark \ref{t-m-tau} we can calculate the values of components
$\{(U_1\ldots U_m \ \zeta^{({m-n},m)})_k\}_{k=0}^{N}$:
\begin{equation}\label{nmUk0}
    (U_1\ldots U_m \ \zeta^{({m-n}, m)})_k =
\end{equation}
\[
=
\begin{cases} \ e^{im\tau\epsilon} \ (g{z})^{m-n-1}\ {g{w} [\alpha_{1}+(g z)^{n} \alpha_2 ]} \   \ & (k=0) \\
e^{im\tau\epsilon} \ [g^2 w^2 (g{z})^{m-n-k-1} \, \alpha_{1} + g^2 w^2 (g{z})^{m-k-1} \, \alpha_2 ] \  \ & (1\leqslant k < m-n)
\\
e^{im\tau\epsilon} \ [g \overline{z} \, \alpha_{1} + g^2 w^2\, (g{z})^{m-k-1}\, \alpha_2 ] \  \ & (k = m-n) \\
e^{im\tau\epsilon} \  g^2 w^2\, (g{z})^{m-k-1}\, \alpha_2 \   \ & (m-n<k<m)  \\
e^{im\tau\epsilon} \ g \overline{z} \, \alpha_{2} \  \ & (k = m) \\
0 \ \ & (m<k\leqslant N)
\end{cases}.
\]
Then, we obtain for (\ref{st-Smn-Sm-t}) :
\begin{eqnarray}\label{Expect-W-alpha-nm}
&&\omega_{\mathcal{S}_{m-n} + \mathcal{S}_{m}}^{m\tau}(\widehat{w}(\alpha_{1})\otimes\widehat{w}(\alpha_2))  \\
&&= \exp\Big[- \frac{1}{4} |w|^2 |\alpha_{1} + (g{z})^{n+1} \alpha_2|^2 |z|^{2(m-n-1)} \frac{1+e^{-\beta_0}}{1-e^{-\beta_0}}\Big]
\nonumber \\
&&\times \exp\Big[- \frac{1}{4} \{|w|^2 (1 - |z|^{2(m-n-1)}) + |z|^2\} |\alpha_{1}|^2
\frac{1+e^{-\beta}}{1-e^{-\beta}}\Big] \nonumber \\
&&\times \exp\Big[- \frac{1}{4} (1 - |w|^2|z|^{2(m-1)}) |\alpha_2|^2  \frac{1+e^{-\beta}}{1-e^{-\beta}}\Big] \ . \nonumber \\
&& {\longrightarrow \exp\Big[- \frac{1}{4} (|\alpha_1|^2 + |\alpha_2|^2 ) \frac{1+e^{-\beta}}{1-e^{-\beta}}\Big]},
\nonumber
\end{eqnarray}
as $ m \to \infty$ {for any fixed $n$}.

Therefore, in this limit, the composed subsystem $\mathcal{S}_{m-n} + \mathcal{S}_m$ evolves
from the initial product of two quasi-free equilibrium states (\ref{st-Smn-Sm}) to the \textit{same} final state, although
for a finite $m$ the evolution (\ref{Expect-W-alpha-nm}) is nontrivial. This again easily understandable taking into account
our analysis of \textit{Subsystem} $\mathcal{S}_m$ and \textit{Subsystem} $\mathcal{S} + \mathcal{S}_m$.

Consider now the case of a {\textit{fixed}} $s := m-n \geqslant 1$. Then the limit in (\ref{Expect-W-alpha-nm}) takes the form:
\begin{eqnarray}\label{Expect-W-alpha-nm-limM}
&& \lim_{m \rightarrow \infty}
\omega_{\mathcal{S}_{s} + \mathcal{S}_{m}}^{m\tau}(\widehat{w}(\alpha_{1})\otimes\widehat{w}(\alpha_2)) = \\
&&= \exp\Big[- \frac{1}{4} |w|^2  |z|^{2(s-1)} |\alpha_{1}|^2
\left\{\frac{1+e^{-\beta_0}}{1-e^{-\beta_0}} - \frac{1+e^{-\beta}}{1-e^{-\beta}}\right\}\Big] \times  \nonumber \\
&&\times \exp\Big[- \frac{1}{4}|\alpha_{1}|^2  \frac{1+e^{-\beta}}{1-e^{-\beta}}\Big]
\exp\Big[- \frac{1}{4}|\alpha_{2}|^2  \frac{1+e^{-\beta}}{1-e^{-\beta}}\Big]  \nonumber \\
&&= \exp\Big[- \frac{1}{4}|\alpha_{1}|^2  \frac{1+e^{-\beta^{**}(s\tau)}}{1-e^{-\beta^{**}(s\tau)}}\Big]
\exp\Big[- \frac{1}{4}|\alpha_{2}|^2  \frac{1+e^{-\beta}}{1-e^{-\beta}}\Big] \ ,\nonumber
\end{eqnarray}
where $\beta^{**}(s\tau)$ verifies equation (\ref{**}). Hence, in this case the limit state (\ref{Expect-W-alpha-nm-limM})
is the product of quasi-free Gibbs states with \textit{different} parameters $\beta^{**}(s\tau)$ and $\beta$. This means that
subsystem $\mathcal{S}_{s}$ keeps a \textit{memory} about perturbation at the moment $t=s\tau$, when the parameter
$\beta^*(s\tau)$ (\ref{*}) of subsystem $\mathcal{S}$ was still different from $\beta$.

Note that (\ref{Expect-W-alpha-nm-limM}) coincides with the
product state  (\ref{st-Smn-Sm}) when $s\rightarrow\infty$.
%%%%%%%%%%%%%%%%%%%%%%%%%%%%%%%%%%%%%%%%%%%%%%%%%%%%%%%%%%%%%%%%%%%%%%%%%%%%%%%%%%%%%%%%%%%%%%%%%%%%%%%%%%%%%%%%%%
\smallskip

\noindent
{\textit{Subsystem} $\mathcal{S}_{\sim n}$. To define $\mathcal{S}_{\sim n}$ } for
$ 0 \leqslant n \leqslant k \leqslant  N$, we divide the total system at the moment
$t=k\tau$ into two subsystems: $\mathcal{S}_{n, k} + \mathcal{C}_{n, k}$. Here
\begin{equation}\label{S-nk}
\mathcal{S}_{n, k} := \mathcal{S} + \mathcal{S}_k  + \mathcal{S}_{k-1} + \cdots + \mathcal{S}_{k-n+1} \ , \
( \mathcal{S}_{0, k} := \mathcal{S} ) \ ,
\end{equation}
whereas
\begin{equation}\label{C-nk}
       \mathcal{C}_{n, k}:= \mathcal{S}_{N} + \cdots + \mathcal{S}_{k+1}
     + \mathcal{S}_{k-n} + \cdots + \mathcal{S}_{1} \ ,
\end{equation}
see definitions in Section \ref{setup}.

{We mean that $\mathcal{S}_{\sim n}$ is an entire ``object" whose entity is $ \mathcal{S}_{n, k} $ at the moment
$ t= k\tau \; $ ( $  k = n, n+1, \cdots, N$ ). As time is running, the elementary subsystems $\mathcal{S}_{k}$ in
$\mathcal{S}_{\sim n}$ are replacing. We study the behaviour of $\mathcal{S}_{\sim n}$  for large $t=k\tau$, i.e., the
$k$-dependence of the ``state" of $ \mathcal{S}_{n, k} $ at $t=k\tau$.}

For any fixed $ t = k\tau$ we can decompose the Hilbert space $\mathscr{H}$ into
a tensor product of two subspaces  $\mathscr{H}_s$ and $ \mathscr{H}_c $ :
\[
            \mathscr{H} = \mathscr{H}_s \otimes \mathscr{H}_c \ .
\]
Here $\mathscr{H}_s$ is the Hilbert space of subsystem (\ref{S-nk})
and $\mathscr{H}_c$ corresponds to subsystem (\ref{C-nk}):
\begin{equation}\label{H-s-H-c}
          \mathscr{H}_s := \mathscr{H}_0\otimes
          \bigotimes_{j=1}^n\mathscr{H}_{k-j+1}, \qquad
           \mathscr{H}_c := \bigotimes_{j=1}^{k-n}\mathscr{H}_j
           \otimes \bigotimes_{j=k+1}^{N}\mathscr{H}_j \ .
\end{equation}

For a density matrix $\varrho$ on $\mathscr{H}$, we introduce the \textit{reduced} density matrix $\varrho_{s}$
on $\mathscr{H}_s$ as a partial trace over $ \mathscr{H}_c $:
\begin{equation}\label{S-varrho}
\varrho_{s} := \Tr_{\mathscr{H}_c}\varrho  \ .
\end{equation}

To avoid a possible confusion caused by the fact that all $\mathscr{H}_j , j= 0,1, \ldots$ are identical to
$\mathscr{F}$ {and by the change of components with time},
we treat the Weyl algebra on the subsystem and the corresponding reduced density matrix of
$ \rho \in \mathfrak{C}_{1}(\mathscr{H}) $ in the following way.
On the Fock space $ \mathcal{F}^{\otimes(n+1)} $ for $n \leqslant N$, we consider the Weyl operator
\begin{equation}\label{W-n-tilde}
           W_n(\zeta) = \exp \Big[i\frac{\langle \zeta, \tilde b, \rangle
      + \langle \tilde b, \zeta \rangle}{\sqrt 2} \Big] ,
\end{equation}
where $\zeta \in \C^{n+1}$, $ \tilde b_0, \cdots, \tilde b_n$ and $\tilde b^*_0, \cdots, \tilde b_n^*$
are the annihilation and the
creation operators in $ \mathcal{F}^{\otimes(n+1)} $ satisfying the corresponding CCR, and
\[
     \langle \zeta, \tilde b \rangle = \sum_{j=0}^n \bar{\zeta_j}
 \tilde b_j, \qquad  \langle \tilde b, \zeta \rangle
     = \sum_{j=0}^n \zeta_j\tilde{b}^*_j .
\]
By $\mathscr{A}(\mathscr{F}^{\otimes(n+1)})$, we denote the $C^*$ algebra generated by the Weyl operators (\ref{W-n-tilde}).
For any subset $J \subset \{1, 2, \cdots, N \} $, we define the operation of taking the partial trace
\[
         R^J : \mathfrak{C}_1(\mathscr{F}^{\otimes(N+1)}) \ni \rho
          \longmapsto R^J\rho \in \mathfrak{C}_1(\mathscr{F}^{\otimes(N+1-|J|)})
\]
by
\[
       \omega_{R^J\rho}\big( W_{N-|J|}(\zeta) \big) = \omega_{\rho}\big( W_N(r_J\zeta) \big) \ .
\]
Here the mapping
\[
          r_J : \C^{N+1-|J|} \ni \zeta \longmapsto r_J \zeta \in \C^{N+1}
\]
is defined by
\[
         (r_J\zeta)_j := \begin{cases} \zeta_0 \quad & ( j=0 )   \\
                                       0            & ( j \in J )  \\
               \zeta_{j-|\{i \in J \, | \, i <j \, \}|} & ( $otherwise$ )
                        \end{cases} \ ,
\]
where $|A|$ denotes the number of elements in the set $A$.

Since all {$\mathscr{H}_1, \mathscr{H}_2, \cdots $ are identical to $\mathscr{F}$, we do not care to} distinguish the spaces
\[
       \bigotimes_{j\in \{0, 1, \cdots, N\} \setminus J}\mathscr{H}_j \qquad
       \mbox{and} \qquad \bigotimes_{j\in \{0, 1, \cdots, N\} \setminus J'}\mathscr{H}_j
\]
for $J \ne J'$ but $ |J| = |J'|$, and consider them as the same space $\mathscr{F}^{\otimes(N+1 -|J|)}$.
Instead, we {pay attention to} distinguishing  projections
\[
     \bigotimes_{j=0}^{N}\mathscr{H}_j \quad \longrightarrow
     \bigotimes_{j\in \{0, 1, \cdots, N\} \setminus J}\mathscr{H}_j
\]
for different subsets $J \subset \{1, 2, \cdots, N \} $ with same $|J|$.

Since we regard $\mathcal{S}_{n,k}$ at time $ t=k\tau$ for
$ k=n, n+1, \cdots $ as the result of the time evolution of a \textit{single} subsystem $\mathcal{S}_{\sim n}$,
we define its state at the moment $ t = k \tau $ by the reduced density matrix
$\{\rho_s(k\tau) \}_{k\geqslant n}$ of this subsystem as follows:
\begin{equation}\label{dyn-S}
    \rho_{s}(k\tau):= R^{\{ 1, \cdots, k-n, k+1, \cdots, N \}}
   \big( \rho(k\tau) \big) = R^{\{ 1, \cdots, k-n, k+1, \cdots, N \}}
      T_{k\tau}(\rho) ,
\end{equation}
see (\ref{Solution-state-t}).
Taking into account Lemma \ref{lem42} and identity $\langle r_J\zeta, r_J\zeta \rangle_{\C^{N+1}} =
\langle \zeta, \zeta \rangle_{\C^{N+1-|J|}}$, one readily obtains the following result.
%%%%%%%%%%%%%%%%%%%%%%%%%%%%%%%%%%%%%%%%%%%%%%%%%%% Lemma %%%%%%%%%%%%%%%%%%%%%%%%%%%%%%%%%%%%%%%%%%%%%%%%%%
\begin{lem}\label{lem51} For the initial density matrix (\ref{Gibbsdensity}),
\[
      \omega_{\rho_{s}(k\tau)}(W_n(\zeta)) =
     \omega_{R^{J_{n,k}}\rho(k\tau)}(W_n(\zeta))
\]
\[
     =  \exp\Big[-\frac{|(U_1 \ldots U_k \, r_{J_{n,k}}\zeta)_0|^2}{4}
         \Big(\frac{1+e^{-\beta_0}}{1-e^{-\beta_0}}
       - \frac{1+e^{-\beta}}{1-e^{-\beta}}\Big)
  - \frac{\langle \zeta, \zeta \rangle}{4}
        \frac{1+e^{-\beta}}{1-e^{-\beta}}\Big]
\]
holds, where $J_{n,k} = \{ 1, 2, \cdots, k-n, k+1, \cdots, N \}$.
\end{lem}
%%%%%%%%%%%%%%%%%%%%%%%%%%%%%%%%%%%%%%%%%%%%%%%%%%%%%%%%%%%%%%%%%%%%%%%%%%%%%%%%%%%%%%%%%%%%%%%%%%%%%%%%%%%%%
\bigskip

To study the limit $k \to \infty$ (and $N\to \infty$ satisfying $k \leqslant N$)
for a fixed $n$, we note that $ (U_1 \ldots U_k \, r_{J_{n,k}} \zeta)_0 \to 0$
follows from (\ref{U0}) and {$|z| <1$ (Hypothesis 5)}. Lemma \ref{lem51} implies that
\begin{equation}\label{S-kn}
        \lim_{k\to\infty}\omega_{\rho_{s}(k\tau)}(W_n(\zeta)) =
        \exp\Big[- \frac{\langle \zeta, \zeta \rangle}{4}
        \frac{1+e^{-\beta}}{1-e^{-\beta}}\Big]
         = \omega_{\rho^{(\beta)}_n}(W_n(\zeta)) \  ,
\end{equation}
where
\begin{equation}\label{Lim-State-S-n}
            \rho^{(\beta)}_n = \exp\big[- \beta\sum_{j=0}^n \tilde b^*_j
    \tilde b_j\big] /Z(\beta )^{n+1}
\end{equation}
by the irreducibility of the CCR algebra
$\mathscr{A}(\mathscr{F}^{\otimes(n+1)})$ and by the definition of $Z(\beta ) = (1-e^{-\beta})^{-1}$.

%\noindent
Therefore, we proved the following statement:
%%%%%%%%%%%%%%%%%%%%%%%%%%%%%%%%%%%%%%%%%%%%%%%%%%%%%%%%% Theorem %%%%%%%%%%%%%%%%%%%%%%%%%%%%%%%%%%%%%%%%%%%%%%%%%%%%%
\begin{thm} \label{thm52}
Let the initial state of the total system $\mathcal{S}+\mathcal{C}$ is
defined by the density matrix (\ref{prodG}): $\rho = \rho(\beta, \beta_0-\beta; e)$.
Then for any fixed $n$, the state $\omega_{\rho_{s}(k\tau)}(\cdot)$
of subsystem $\mathcal{S}_{n, k}$ converges to the equilibrium Gibbs state
$\omega_{\rho^{(\beta)}_n}(\cdot)$ {as $k\to\infty$ in the weak*-topology
for the states on  $\mathscr{A}(\mathscr{F}^{\otimes(n+1)})$}, see (\ref{w*-lim}).
\end{thm}
%%%%%%%%%%%%%%%%%%%%%%%%%%%%%%%%%%%%%%%%%%%%%%%%%%%%%%%%%%
%%%%%%%%%%%%%%%%%%%%% Theorem %%%%%%%%%%%%%%%%%%%%%%%%%%
%%%%%%%%%%%%%%%%%%%%%%%%%%%%%%%%%%%%%%
\begin{thm}\label{thm53}
Under the same conditions as in Theorem \ref{thm52}, one gets
\[
        \lim_{k\to\infty}S(\rho_s (k\tau)) = S(\rho^{(\beta)}_n) \ .
\]

\end{thm}
%%%%%%%%%%%%%%%%%%%%%%%%%%%%%%%%%%%%%%%%%%%%%%%%%%%%%%%%%%%%%%%%%%%%%%%%%%%%%%%%%%%%%%%%%%%%%%%%%%%%%%%%%%%%%%%%%%%%%%%%
{\sl Proof }:
Let vector $\xi_{n,k} \in \C^{n+1}$ be defined by
$(U_1\ldots U_k r_{J_{n,k}}\zeta)_0 =: \langle \xi_{n,k}, \zeta \rangle$.
Then $ \langle \xi_{n,k}, \xi_{n,k} \rangle \to 0 $ as $k \to \infty$ for fixed $n$.
By Lemma \ref{lem-W-S} and Lemma \ref{lem51},  we obtain
\[
      S(\rho_s(k\tau)) = n \sigma\Big(
       \frac{1+e^{-\beta}}{1-e^{-\beta}} \Big) +
       \sigma\Big(\frac{1+e^{-\beta}}{1-e^{-\beta}}
        + \langle \xi_{n,k}, \xi_{n,k} \rangle
     \Big(\frac{1+e^{-\beta_0}}{1-e^{-\beta_0}}
           - \frac{1+e^{-\beta}}{1-e^{-\beta}}\Big)\Big)
\]
\[
          \longrightarrow (n+1)\sigma\Big(
       \frac{1+e^{-\beta}}{1-e^{-\beta}} \Big) = S(\rho_n^{(\beta)}).
\]
\hfill$\square$

\begin{remark} The local entropy decreases or increases  according
to $\beta > \beta_0$ or $\beta < \beta_0$, respectively.
\end{remark}

%\newpage
%%%%%%%%%%%%%%%%%%%%%%%%%%%%%%%%%%%%%%%%%% subsection Short-time-Limit %%%%%%%%%%%%%%%%%%%%%%%%%%%%%%%%%%%%%%%%%%%%%%%%
\subsection{A Short-Time Limit for Repeated Perturbations}\label{STL}
%%%%%%%%%%%%%%%%%%%%%%%%%%%%%%%%%%%%%%%%%%%%%%%%%%%%%%%%%%%%%%%%%%%%%%%%%%%%%%%%%%%%%%%%%%%%%%%%%%%%%%%%%%%%%%%%%%%%%%%
The results in the previous Section \ref{subsystem} are essentially due our explicit knowledge of the initial
density matrix (\ref{state-S-C}), (\ref{Gibbsdensity}) of the total system $\mathcal{S} + \mathcal{C}$.
In this subsection, we show that the lack of this information is not decisive for certain results concerning the convergence to
equilibrium if one considers a short-time limit for the repeated interactions.

Let us study it for example of the subsystem $\mathcal{S}$.  We keep to consider the initial state
of the system $\mathcal{S} + \mathcal{C}$ to be a product state with the density matrix
\begin{equation}\label{0-product-state}
\rho = \rho_{0}\otimes \bigotimes_{k=1}^N \rho_k  \in  \mathfrak{C}_{1}(\mathscr{H}) \ ,
\end{equation}
see (\ref{state-S-C}), but we essentially relax the conditions on $\rho_{0}$ and $\{\rho_k\}_{k=1}^N$:
\begin{eqnarray*}
\mbox{(h1)}& \qquad    \rho_1 =\rho_2 = \cdots =\rho_N  \in  \mathfrak{C}_{1}(\mathscr{F}) \ ; \\
\mbox{(h2)}&
       { \Tr_{\mathscr{F}} [\rho_1 a] = \Tr_{\mathscr{F}}  [\rho_1 a^2] =
        \Tr_{\mathscr{F}}  [\rho_1 a^*]
= \Tr_{\mathscr{F}}  [\rho_1 a^{*2}] = 0 }  ; \\
\mbox{(h3)}& \qquad
{\Tr_{\mathscr{F}}  [\rho_1 (a^*a)^2]} < \infty  \ .
\end{eqnarray*}
%%%%%%%%%%%%%%%%%%%%%%%%%%%%%%%%%%%%%%%%%%%%%%%%%%%% Remark %%%%%%%%%%%%%%%%%%%%%%%%%%%%%%%%%%%%%%%%%%%%%%
\begin{remark}\label{S-TRPLimGibbs}
{Note that hypothesis (h1)-(h{3}) are satisfied when the density matrices $\{\rho_k\}_{k=0}^N$ correspond to the
gauge-invariant quasi-free states with parameter $\beta_0$ for $k=0$ and $\beta$ for $k=1,2,\ldots,N$, see
(\ref{Gibbsdensity}). Then (h2) is due to the gauge invariance and one gets for (h3):
\begin{equation}\label{major}
\Tr_{\mathscr{F}}[\rho_k (a^*a)^2] = (2n_{\beta}^2+ n_{\beta}) \ ,
\end{equation}
{where} $n_{\beta}= \Tr_{\mathscr{F}}\rho_k (a^*a)=(e^{\beta}-1)^{-1}$, $k=1, \ldots, N$.}
\end{remark}
%%%%%%%%%%%%%%%%%%%%%%%%%%%%%%%%%%%%%%%%%%%%%%%%%%%%%%%%%%%%%%%%%%%%%%%%%%%%%%%%%%%%%%%%%%%%%%%%%%%%%%%%%%

Below we denote by $|ya^* + \bar ya|$ the operator originated from the \textit{polar decomposition} of the self-adjoint
operator {$ya^* + \bar ya = U \, |ya^* + \bar ya| $}, where $U$ is the partial isometry {on} $\mathscr{F}$.
%%%%%%%%%%%%%%%%%%%%%%%%%%%%%%%%%%%%%%%%%%%%%%% Lemma %%%%%%%%%%%%%%%%%%%%%%%%%%%%%%%%%%%%%%%%%%%%%%%%%%%%%%
\begin{lem} \label{lem61}
Under hypothesis (h1)-(h3), the following bounds hold:
\begin{eqnarray*}
\rm{(i)} \hspace{20mm}&
                 \Tr_{\mathscr{F}} [\rho_k a^*a ] < \infty  ,
\\
\rm{(ii)} \hspace{20mm} &
                  \Tr_{\mathscr{F}} [\rho_k |ya^* + \bar ya|^2] \leqslant C|y|^2 ,
\\
({\rm iii}) \hspace{20mm} &
                    \Tr_{\mathscr{F}} [\rho_k |ya^* + \bar ya|^3] \leqslant C'|y|^3 ,
\\
({\rm iv}) \hspace{20mm} &
                    \Tr_{\mathscr{F}} [\rho_k |ya^* + \bar ya|^4] \leqslant C''|y|^4 ,
\end{eqnarray*}
for all {$k = 1, \ldots , N$}. Here  $ C, C', C''$ are positive constants, which
depend only on {$\Tr[ \rho_1(a^*a)^2]$}.
\end{lem}
%%%%%%%%%%%%%%%%%%%%%%%%%%%%%%%%%%%%%%%%%%%%%%%%%%%%%%%%%%%%%%%%%%%%%%%%%%%%%%%%%%%%%%%%%%%%%
{\sl Proof : } The first bound (i) is a consequence of the Cauchy-Schwarz
inequality and (h3). Applying the inequalities
\[
       |A + A^*|^2 \leqslant |A + A^*|^2 + |A - A^*|^2 = 2(AA^* + A^*A),
\]
\[
        |A + A^*|^4 \leqslant |A + A^*|^4 + |A - A^*|^4
             + |A + iA^*|^4 + |A -i A^*|^4
\]
\[
       = 4(AA^* + A^*A)^2 + 4(A^2A^{*2} + A^{*2}A^2),
\]
to $ A= \bar y a$, we obtain (ii) and (iv). Finally, a combination of (ii), (iv) with the Cauchy-Schwarz inequality
yields (iii). \hfill $\square$
%\medskip
%%%%%%%%%%%%%%%%%%%%%%%%%%%%%%%%%%%%%%%%%%%%%%%%%% Theorem %%%%%%%%%%%%%%%%%%%%%%%%%%%%%%%%%%%%%%%%%%%%%%%%%%
\begin{thm}\label{thm62}
Let $\tau \rightarrow 0$, $N \to \infty$ be short-time perturbation limit {subject} to demands:
$\tau^2N\to \infty$ and $\tau^3N \to 0$.
Then for any initial condition (\ref{0-product-state}) verifying (h1)-(h3),
{ the characteristic function $\omega_{\mathcal{S}}^{N\tau}(\widehat{w}(\theta))$ of the state for subsystem
$\mathcal{S}$ at $t=N\tau$, converges to}:
\begin{eqnarray}\label{tau-N-lim}
&&\omega_{\mathcal{S}}(\widehat{w}(\theta)): =
\lim_{\tau \rightarrow 0, N \to \infty}\omega_{\mathcal{S}}^{N\tau}(\widehat{w}(\theta)) =\\
&&\lim_{\tau \rightarrow 0, N \to \infty} \omega_{\rho(N\tau)}(W(\zeta_{\theta}))=
e^{- |\theta|^2 \Tr_{\mathscr{F}}[\rho_{1}\,(a^*a +aa^*)]/4}
\nonumber \ .
\end{eqnarray}
Here $\theta \in \C $ and the $(N+1)$-component vector for the $\mathcal{S} + \mathcal{C}$ Weyl operator is
\[
           \zeta_{\theta} = \begin{pmatrix} \theta \\
         0 \\ \cdot \\ \cdot \\
               \cdot \\ 0
                      \end{pmatrix} \in \C^{N+1} \ .
\]
\end{thm}
%%%%%%%%%%%%%%%%%%%%%%%%%%%%%%%%%%%%%%%%%%%%%%%%%%%% Remark %%%%%%%%%%%%%%%%%%%%%%%%%%%%%%%%%%%%%%%%%%%%%%
\begin{remark}\label{thm62R}
By (\ref{tau-N-lim}) the state $\omega_{\mathcal{S}}^{N\tau}$ converges to $\omega_{\mathcal{S}}$ in the weak*-topology,
see Appendix \ref{App1}, A.4. From the right-hand side of (\ref{tau-N-lim}) and Definition \ref{QF-state} we deduce that
the limit state is gauge-invariant and quasi-free with $h(\theta):= |\theta|^{2} \, \Tr_{\mathscr{F}}[\rho_{1} \, a^*a]$.
\end{remark}
%%%%%%%%%%%%%%%%%%%%%%%%%%%%%%%%%%%%%%%%%%%%%%%%%%%% Remark %%%%%%%%%%%%%%%%%%%%%%%%%%%%%%%%%%%%%%%%%%%%%%
\begin{remark}\label{thm62R-T}
Recall that the state $\omega$ over the Weyl algebra $\mathscr{A}(\mathscr{F})=\overline{\mathscr{A}_{w}(\mathscr{F})}$
is \textit{regular}, $C^n$-smooth or analytic, if the function (see (\ref{S-Weyl}))
\begin{equation}\label{s-function}
s \mapsto \omega (\widehat{w}( s \theta)) = \omega (e^{i \, s {\Phi(\theta)}/\sqrt 2})
\end{equation}
is respectively continuous, $C^n$-smooth or analytic in the vicinity of $s = 0$.
In the {last} case the characteristic function
$\omega (\widehat{w}(s \theta))$ (and therefore the state) is completely determined by
\begin{equation}\label{T-repr}
\omega (\widehat{w}(s \theta)) =
\exp \left\{\sum_{m=1}^{\infty} \frac{i^m s^m}{m!} \, 2^{-m/2} \ \omega^{T}(\Phi^m (\theta))\right\} \ .
\end{equation}
Here $\{\omega^{T}(\Phi^{m}(\theta))\}_{m=0}^{\infty}$ are \textit{truncated} correlation functions defined
recursively by relations \cite{BR2}, \cite{Ve}:
\begin{eqnarray*}
%\label{truncated}
&& \omega^{T}(\Phi(\theta)): = \omega(\Phi(\theta)) \ , \\
&& \omega^{T}(\Phi^2(\theta)) := \omega(\Phi^2(\theta)) - \omega(\Phi(\theta))^2 \ ,  \nonumber \\
&& \omega^{T}(\Phi^3(\theta)) := \omega(\Phi^3(\theta)) - 3 \omega(\Phi^2(\theta)) \omega(\Phi(\theta)) +
2 \omega(\Phi(\theta))^3 \, , \ {\rm{etc}} \nonumber
\end{eqnarray*}
Lemma \ref{lem61} implies that states corresponding to density matrices {$\rho_1= \rho_2 = \ldots$} are $C^4$-smooth.
\end{remark}
%%%%%%%%%%%%%%%%%%%%%%%%%%%%%%%%%%%%%%%%%%%%%%%%%%%%%%%%%%%%%%%%%%%%%%%%%%%%%%%%%%%%%%%%%%%%%%%%%%%%%%%%%%%%%%%%%%%%
{\sl Proof } (of Theorem \ref{thm62}): By (h2) and by Lemma \ref{lem61} (i)-(iii) together with  Remark \ref{thm62R-T},
we obtain for the states $\omega(\cdot) = \omega_{\rho_k}(\cdot)$ the representation of (\ref{T-repr}) in the form:
\begin{eqnarray} \label{Ck}
C_k(\theta) = \omega_{\rho_k} (\widehat{w}(\theta)) =  \exp [-\frac{1}{4} \ \omega^{T}_{\rho_k}(\Phi^2(\theta)) + R(\theta)] \ , \
k = 1,2, \ldots , N ,
\end{eqnarray}
where $R(\theta)= O(|\theta|^3)$ in the vicinity of $\theta = 0$. {For} the self-adjoint operator
$\Phi(\theta) = \bar \theta a + \theta a^*$, the hypothesis (h2) and Lemma \ref{lem61} (i) imply
\begin{equation}\label{omega-T}
\omega^{T}_{\rho_k}(\Phi^2(\theta)) = |\theta|^2 \ \Tr_{\mathscr{F}}[\rho_k \, (a^*a+aa^*)] \ .
\end{equation}

Now, taking into account Lemma \ref{lem21} for the vector $\zeta_{\theta}$, (\ref{Ck}) {and}
(\ref{omega-T}),  we obtain the representation:
\[
      {\omega_{\mathcal{S}}^{N\tau}(\widehat{w}(\theta))} =
       \omega_{\rho(N\tau)}(W(\zeta_{\theta})) =
       C_0(e^{i\epsilon\tau N}(gz)^N\theta)\prod_{k=1}^N
       C_k({e^{i\epsilon\tau N}g  w }\ (gz)^{N-k}\theta)
\]
\begin{equation}\label{tau-N}
 = C_0(e^{i\epsilon\tau N}(gz)^N\theta) \exp\Big(-\sum_{k=1}^N\frac{|\theta_k|^2}{4}
\Tr_{\mathscr{F}}[(a^*a+aa^*)\rho_k]  + \widehat{R} \Big) \ .
\end{equation}
Here by (\ref{Uk}) and  by (\ref{Ck}) one has
\[
\theta_k := {e^{i\epsilon N\tau}g w} \ (gz)^{N-k}\theta \ , \
\sum_{k=1}^N|\theta_k|^2 = |\theta|^2 |w|^2 \frac{1-|z|^{2N}}{1-|z|^2} \ ,\ \widehat{R} = \sum_{k=1}^N O(|\theta_k|^3) \ .
\]

{By virtue of} (\ref{g-w}) and (\ref{z}), we get $ |g(\tau)| = 1$, $|w(\tau)|^2 + |z(\tau)|^2 =1$  and  also
\[
w(\tau) = i\eta \tau +  O(\tau^3)\ , \ |z(\tau)| = 1 - \frac{|\eta|^2\tau^2}{2} + O(\tau^4)  \ ,
\]
for small $\tau$.
This yields for {small $\tau > 0$ and large $N$}, the estimates
$|(gz)^N| \leq  O(e^{-|\eta|^2\tau^2N/2})$, $|\theta_k| \leq  O(\tau)$, and  $\widehat{R}=O(\tau^3 N)$ by virtue of (h1).
Then taking into account the conditions $\tau^2N\to \infty$ and $\tau^3N \to 0$, we get the limits:
\[
\lim_{\tau \rightarrow 0, N \to \infty} C_0(e^{i\epsilon\tau N}(gz)^N\theta) = 1 \ , \
\lim_{\tau \rightarrow 0, N \to \infty} \sum_{k=1}^N|\theta_k|^2 = |\theta|^2 \ , \
\lim_{\tau \rightarrow 0, N \to \infty}  \widehat{R} = 0 \ .
\]
Note that $C_0$ is a continuous function  because it {is} defined by a normal state with density matrix $\rho_0$,
see {(\ref{C-k})}.

Inserting all these limits into (\ref{tau-N}), we obtain what is claimed as the limit (\ref{tau-N-lim}).    \hfill $\square$

%%%%%%%%%%%%%%%%%%%%%%%%%%%%%%%%%%%%%%%%%%%%%%%%%%%% Corollary %%%%%%%%%%%%%%%%%%%%%%%%%%%%%%%%%%%%%%%%%%%%%%%%%%%%%%%
\begin{cor}\label{thm62cor}
Suppose (see Section \ref{RepPertQFDyn}) that all $\{\rho_k\}_{k=1}^N$ correspond to the gauge-invariant
quasi-free Gibbs state with parameter $\beta$
(\ref{Gibbsdensity}):
\[
           \rho_k = e^{-\beta a^*a}/\Tr_{\mathscr{F}}[e^{-\beta a^*a}]
       \ , \quad (k = 1,2,\ldots, N ) \ .
\]
These states satisfy (h1)-(h3).
The statement in Theorem \ref{thm62} is valid with the limit
\begin{eqnarray}\label{tau-N-lim-beta}
\lim_{\tau \rightarrow 0, N \to \infty}
\omega_{\mathcal{S}}^{N\tau}(\widehat{w}(\theta))=
\exp \left\{- \frac{|\theta|^2}{4} \frac{1+e^{-\beta}}{1 - e^{-\beta}}\right\}  \, .
\end{eqnarray}
It coincides with the result for equilibrium state (\ref{st-S-t=inf}) of the subsystem $\mathcal{S}$
when the finite step $\tau$ verifies the Hypothesis 5.
\end{cor}
%%%%%%%%%%%%%%%%%%%%%%%%%%%%%%%%%%%%%%%%%%%%%%%%%%%%%%%%%%%%%%%%%%%%%%%%%%%%%%%%%%%%%%%%%%%%%%%%%%%%%%%%%%%%%%%%%%%%%%

{Note that our choice of the short-time perturbation limit $\tau \rightarrow 0$, $N \to \infty$  subjected to
$\tau^2N\to \infty$ and $\tau^3N \to 0$ gives a \textit{universal} gauge-invariant quasi-free limiting state
under hypothesis (h1)-(h{3}). The hypotheses (h2), (h3) control only first "two moments" in the creation-annihilations
operators of the reference initial states of subsystem $\mathcal{C}$. Together with stationarity and independence of repeated
perturbations in (h1), these conditions make our observation similar to well-known universal laws similar to the
non-commutative Central Limit Theorem \cite{Ve}.}

This similarity is bolstered by the fact that the state $\omega_{\rho_0}$ of the subsystem $\mathcal{S}$ may be replaced
by any {\textit{regular} state}. This indicates how large {could be} the ``basin of attraction" of the universal limiting
gauge-invariant quasi-free state. Another common point is the method of characteristic functions relevant in the both
cases \cite{Ve}.

%%%%%%%%%%%%%%%%%%%%%%%%%%%%%%%%%%%%%%%%%%%%%%%%% Acknowledgments %%%%%%%%%%%%%%%%%%%%%%%%%%%%%%%%%%%%%%%%%%%%%%%%
\vspace{1cm}

\noindent
\textbf{Acknowledgements }

\noindent {H.T. thanks  JSPS for the financial support under the Grant-in-Aid for Scientific Research (C) 24540168.
He is also grateful to Aix-Marseille and Toulon Universities for their hospitality.}

V.A.Z. is thankful to  Marco Merkli for instructive conversations and remarks.

%%%%%%%%%%%%%%%%%%%%%%%%%%%%%%%%%%%%%%%%%%%%%%%%%%%%%%%%%%%%%%%%%%%%%%%%%%%%%%%%%%%%%%%%%%%%%%%%%%%%%%%%%%%%%%%%%%%
%%%%%%%%%%%%%%%%%%%%%%%%%%%%%%%%%%%%%%%%%%%%%%%%% Appendix %%%%%%%%%%%%%%%%%%%%%%%%%%%%%%%%%%%%%%%%%%%%%%%%%%%%%%%%
%\newpage
%%%%%%%%%%%%%%%%%%%%%%%%%%%%%%%%%%%%%%%%%%%%%%%%%%%%%%%%%%%%%%%%%%%%%%%
%%%%%%%%%%%%%%%%%%%%%%%%%%%%%%%%%%%%%%%%%%%%%%%
\appendix
\section{Appendix}\label{App1}
%%%%%%%%%%%%%%%%%%%%%%%%%%%%%%%%%%%%%%%%%%%%%%%
%%%%%%%%%%%%%%%%%%%%%%%%%%%%%%%%%%%%%%%%%%%%%%%%%%%%%%%%%%%%%%%%%%%%%%%%%%%%%%%%%%%%%%%%%%%%%%%%%%%%%%%%%%%%%%%%%%%%%
Let $\omega_{\varrho(t)}(\cdot)$ be time dependent normal state $\Tr_{\mathscr{H}} (T_{t}(\varrho) \ \cdot \ )$ on
$\mathcal{L}(\mathscr{H})$. Then evolution operator $T_{t}: \varrho \mapsto U(t)\varrho\, U^*(t)$ for $t\in \mathbb{R}$ on
the set of density matrices (\ref{Solution-state-t}) defines a dual $\ast$-automorphisms $T_{t}^*$ of the $C^*$-algebra of bounded
operator {$\mathcal{L}(\mathscr{H})$} (the Heisenberg picture on the $W^*$-algebra). We collect here some general remarks about
the $C^*$-dynamical systems versus $W^*$-dynamical setting, cf Remark \ref{rem-dual}.
%%%%%%%%%%%%%%%%%
\smallskip

\noindent{\bf A1.} Recall that the \textit{dual} of the Banach space $\mathfrak{C}_{1}(\mathscr{H})$ is the set
$\mathfrak{C}_{1}(\mathscr{H})^*:=\mathcal{L}(\mathfrak{C}_{1}(\mathscr{H}), \mathbb{C})$ of \textit{all} bounded linear
functionals on $\mathfrak{C}_{1}(\mathscr{H})$, and that they coincide with
$\{\phi \mapsto \Tr_{\mathscr{H}} (\phi \, A )\}_{A \in \mathcal{L}(\mathscr{H})}$. Here the {correspondence} is an isometric
isomorphism of $\mathcal{L}(\mathscr{H})$ onto $\mathfrak{C}_{1}(\mathscr{H})^*$ such that the norm of each function is
${\|\Tr_{\mathscr{H}} (\, \cdot \, A )\|_{\mathfrak{C}_{1}^{*}}} = \|A\|$, or
$\mathcal{L}(\mathscr{H})=\mathfrak{C}_{1}(\mathscr{H})^*$.

Note that the set of maps $\{A \mapsto \Tr_{\mathscr{H}} (\phi \, A )\}_{\phi \in \mathfrak{C}_{1}(\mathscr{H})}$
{does not cover the set of all} continuous linear functionals on $\mathcal{L}(\mathscr{H})$,
but they yield \textit{dual} of the subspace, {which consists} of the compact operators:
$\mathfrak{C}_{\infty}(\mathscr{H})^{*}=\mathfrak{C}_{1}(\mathscr{H})$, with the norm of each functional:
${\|\Tr_{\mathscr{H}} (\phi \; \cdot \, )\|_{\mathfrak{C}_{\infty}^{*}}} = \|\phi\|_1$.

Therefore, to control the $\|\cdot\|_{1}$-continuity of \textit{density matrix} evolution $T_{t}(\varrho)$
{with help of duality (\ref{dual})}, one needs some additional arguments. To this end, note that since dynamics
(\ref{Solution-state-t}) is trace- and positivity- preserving, one gets $\|T_{t}(\varrho)\|_{1} =1$.
By {the strong continuity of (\ref{U})} (or by the unity-preserving $T_{t}^*(\mathbb{1}) = \mathbb{1}$, and duality)
we also get the \textit{weak operator} continuity of $T_{t}(\varrho)$. Together these arguments yield the
$\|\cdot\|_{1}$-continuity of $T_{t}(\varrho)$ , see e.g. \cite{Za}, Ch.2.4.

\smallskip

\noindent{\bf A2.} \textit{$C^*$-dynamical systems.} It is a pair $(\mathfrak{A} , \tau^t )$, where $\mathfrak{A}$ is a unital
$C^*$-algebra and $\tau^t $ is a strongly continuous (continuous in topology of this algebra) group of $\ast$-automorphisms
of $\mathfrak{A}$, see \cite{AJP1}, \cite{BR1}.

In the context of Remark \ref{rem-dual} one identifies $\mathfrak{A}$ with $\mathscr{A(H)}$ (or $\mathcal{L}(\mathscr{H})$)
and $\tau^t$  with dynamics $T_{t}^*$ dual with respect to the state
$\omega_{\varrho}(\cdot):= \Tr_{\mathscr{H}} (\varrho \ \cdot )$ on $\mathcal{L}(\mathscr{H})$.
Here {$\varrho$ is a density matrix} (see A1) and
\begin{equation}\label{def-T}
\omega_{\varrho}(T_{t}^*(A)) = \omega_{\, T_{t}(\varrho)}(A) \ .
\end{equation}

In the case of boson systems (Section \ref{HDIS}), the $C^*$-approach is too restrictive.
First, it is because the CCR force us to use the Weyl algebra $\mathscr{A}(\mathscr{H})$ and this $C^*$-algebra is only a
\textit{subalgebra} of $\mathcal{L}(\mathscr{H})$. Hence, a \textit{predual} to $\mathscr{A}(\mathscr{H})$ is not
$\mathfrak{C}_{1}(\mathscr{H})$, see A1.
Second, since the CCR beak the operator-norm continuity of dynamics $T_{t}^*(W(\zeta))= W(U(t)\zeta)$:
\begin{equation*}\label{}
{\| W(\zeta_1) - W(\zeta_2)\|_{\mathscr{A}(\mathscr{H})} = 2 \ ,
\ {\rm{if}} \  \zeta_1 \neq \zeta_2 \ }.
\end{equation*}

%%%%%%%%%%%%%%%%%%%%%%%%%%%%%%%%%%%%% Remark %%%%%%%%%%%%%%%%%%%%%%%%%%%%%%%%%%%%%
%\begin{remark}\label{W*}
\smallskip

\noindent{\bf A3.} \textit{$W^*$-dynamical systems.} To overcome difficulties mentioned in A2, one has to take a closure,
{$\mathfrak{M}(\mathscr{H})$}, of the Weyl algebra {$\mathscr{A}(\mathscr{H})$} (\ref{S-Weyl}) in {topology} which is weaker
than the operator-norm topology of this $C^*$-algebra.

To this end, consider on $\mathcal{L}(\mathscr{H})$ the weak*-topology ($w^\ast$-topology) generated by the set of
linear functionals $\{ A \mapsto \Tr_{\mathscr{H}} (\phi \, A )\}_{\phi \in \mathfrak{C}_{1}(\mathscr{H})}$,
 see A1.
A priori, it is stronger than the weak operator topology on $\mathcal{L}(\mathscr{H})$, but weaker
than the operator-norm or the weak Banach space topology on $\mathcal{L}(\mathscr{H})$ since
$\mathfrak{C}_{1}(\mathscr{H})\subset \mathcal{L}(\mathscr{H})^\ast$.
By A1 the trace-class $\mathfrak{C}_{1}(\mathscr{H})= \mathcal{L}(\mathscr{H})_{\ast}$ is \textit{predual} of
$\mathcal{L}(\mathscr{H})$ since $\mathfrak{C}_{1}(\mathscr{H})^*= \mathcal{L}(\mathscr{H})$.

If we denote by {$\mathfrak{M}(\mathscr{H})$} the closure of the Weyl algebra
{$\mathscr{A}(\mathscr{H})$} (\ref{Wz-bis}) in the $w^\ast$-topology, then it is the von Neumann algebra acting on
the boson Fock space $\mathscr{H}$.
Note that $\mathfrak{M}(\mathscr{H})$ {is} $\ast$-isomorphic to  $\mathcal{L}(\mathscr{H})$.
By construction of the von Neumann algebra and
by duality (\ref{def-T}) the $\ast$-automorphism $t \mapsto T_{t}^*(A)$ of
{$\mathfrak{M}(\mathscr{H})$} is continuous in the $w^\ast$-topology ($W^*$-dynamics).

A $W^*$-dynamical system is a pair $(\mathfrak{M}, T_{t}^*)$, where $\mathfrak{M}$ is a von Neumann algebra acting on
a Hilbert space $\mathscr{H}$ and $T_{t}^*$ is a  $W^*$-dynamics on $\mathfrak{M}$, see  e.g. \cite{AJP1}, \cite{BR1}
for details.

%%%%%%%%%%%%%%%%%%%%%%%%%%%%%%%%%%%%%%%%%%%%%% Remark %%%%%%%%%%%%%%%%%%%%%%%%%%%%%%%%%%%%%%%%%%%
\smallskip

\noindent{\bf A4.} Note that the finite linear combinations of (\ref{Wz}) are norm dense in the Weyl $C^*$-algebra $\mathscr{A(H)}$
and the map: $\zeta \mapsto W(\zeta)$ is continuous in the strong operator topology. Then by the Araki-Segal theorem,
see e.g. \cite{AJP1}, any state $\omega$ on this algebra is completely determined by its characteristic function
\begin{equation}\label{E-W}
\zeta \mapsto E_{\omega}(\zeta): = \omega (W(\zeta)) \ .
\end{equation}
{When} the function $s \mapsto E_{\omega}(s \, \zeta)$ is continuous, the state $\omega$ is \textit{regular}.
Recall that the smoothness of this function near $s = 0$ decides the $C^m$-smoothness or analyticity of $\omega$, see
Remark \ref{thm62R-T}.

The set of states $S_{\mathcal{A}}$ over algebra $\mathscr{A(H)}$ is a subset of a dual to this
algebra: $S_{\mathcal{A}} \subset \mathscr{A(H)}^{\ast}$.
Besides the uniform topology on $\mathscr{A(H)}^{\ast}$ one considers also the weak*-topology.
Restriction of this topology to $S_{\mathcal{A}}$ is  defined by the base of neighbourhoods
\begin{equation}\label{weak-star-top}
\mathcal{N}(\omega; A_1,\ldots, A_n):= \{\omega^{\prime}, \in \mathscr{A(H)}^{\ast}:
|\omega^{\prime}(A_i) - \omega (A_i)| < \varepsilon , \; i = 1, 2, \ldots, n\}
\end{equation}
for any $\varepsilon > 0$ and finite sets of operators $A_1, A_2, \ldots, A_n \in \mathscr{A(H)}$.
If the sequence of {regular} states $\{\omega^{(k)}\}_{k\geq 1}$ {enjoys} the convergence
of characteristic functions
\begin{equation}\label{lim-E-k}
\lim_{k \rightarrow \infty} E_{\omega^{(k)}}(\zeta) = E_{\infty}(\zeta) \ , \ \zeta \in \mathbb{C} \ ,
\end{equation}
then $E_{\infty}(\zeta)$ verifies conditions of the Araki-Segal theorem and defines on $\mathscr{A(H)}$ a regular
state $\omega^{(\infty)}$: $E_{\infty}(\zeta) = E_{\omega^{(\infty)}}(\zeta)$.
By definition (\ref{E-W}) and by (\ref{lim-E-k})
this state is the limit of the sequence $\{\omega^{(k)}\}_{k\geq 1}$ in the weak*-topology (\ref{weak-star-top}) on
$S_{\mathcal{A}}$:
\begin{equation}\label{w*-lim}
\omega^{(\infty)} = w^{\ast} \!-\! \lim_{k \rightarrow \infty}  \omega^{(k)} \ ,
\end{equation}
see e.g. \cite{AJP1} and \cite{BR2}, Ch.5.2.5.

%%%%%%%%%%%%%%%%%%%%%%%%%%%%%%%%%%%%%%%%%%%%%%%%% Refererences %%%%%%%%%%%%%%%%%%%%%%%%%%%%%%%%%%
\newpage
%%%%%%%%%%%%%%%%%%%%%%%%%%%%%%%%%%%%%%%%%%%%%%%%%%%%%%%%%%%%%%%%%%%%%%%

%%%%%%%%%%%%%%%%%%%%%%%%%%%%%%%%%%%%%%%%%%%%%%%%%%%%%%%%%%%%%%%%%%%%%%%%%%%%%%%%%%%%%%%%%%%%%%%%%%

%%%%%%%%%%%%%%%%%%%%%%%%%%%%%%%%%%%%%%%%%%%%%%%%%%%%%%%%%%%%%%%%%%%%%%%%%%%%%%%%%%%%%%%%%%%%%%%%%%%%%%%%%%%%%%%%%%%%

\begin{thebibliography}{999999}
%%%%%%%%%%%%%%%%%%%%%%%%%%%%%%%%%%%%%%%%%%%%%%%%%%%%%%%%%%%%%%%%%%%%%%%

\bibitem[Ar0]{Ar0} H.Araki, On the Diagonalization of a Bilinear Hamiltonian by a Bogoliubov Transformation,
\textit{Publ. RIMS, Kyoto Univ. Ser.A}, \textbf{4} (1968), 387-412.
\bibitem[Ar1]{Ar1} H.Araki, Relative Entropy of States of von Neumann Algebras,
\textit{Publ. RIMS, Kyoto Univ.}, \textbf{11} (1976), 809-833.
\bibitem[AJP1]{AJP1} \textit{Open Quantum Systems I, The Hamiltonian Approach},
S.~Attal, A.~Joye, C.-A.~Pillet (Eds.),
Lecture Notes in Mathematics \textbf{1880}, Springer-Verlag, Berlin-Heidelberg 2006.
\bibitem[AJPII]{AJP2} \textit{Open Quantum Systems II, The Markovian Approach},
S.~Attal, A.~Joye, C.-A.~Pillet (Eds.),
Lecture Notes in Mathematics \textbf{1881}, Springer-Verlag, Berlin-Heidelberg 2006.
\bibitem[AJP3]{AJP3} \textit{Open Quantum Systems III, Recent Developements},
S.~Attal, A.~Joye, C.-A.~Pillet (Eds.),
Lecture Notes in Mathematics \textbf{1882}, Springer-Verlag, Berlin-Heidelberg 2006.
\bibitem[BJM]{BJM} L.Bruneau, A.Joye, and M.Merkli, Repeated interactions in open quantum systems,
(May 14, 2013). Submitted to \textit{J.Math.Phys. }
\bibitem[BR1]{BR1} O.~Bratteli and D.W.~Robinson,
\textit{Operator Algebras and Quantum Statistical Mechanics}, vol.1, Springer-Verlag, Berlin 1979.
\bibitem[BR2]{BR2} O.~Bratteli and D.W.~Robinson,
\textit{Operator Algebras and Quantum Statistical Mechanics},vol.2, Springer-Verlag (2nd Edt), Berlin 1997.
\bibitem[Fa]{Fa} M.Fannes, The entropy of quasi-free states for a continuous boson system, \textit{Ann.de l'IHP}, section A,
\textbf{28}(1978) 187-196.
\bibitem[NVZ]{NVZ} B. Nachtergaele, A. Vershynina, and V. A. Zagrebnov, \\
Non-Equilibrium States of a Photon Cavity Pumped by an Atomic Beam,
\textit{Annales Henri Poincar\'{e}}, \textbf{15} (2014), 213-262.
\bibitem[Ve]{Ve} A.F.~Verbeure,
\textit{Many-Body Boson Systems}, Springer-Verlag, Berlin 2011.
\bibitem[Za]{Za} V.A.~Zagrebnov,
\textit{Topics in the Theory of Gibbs Semigroups}, KU Leuven University Press, Leuven 2003.
\end{thebibliography}
\end{document}